\documentclass[10pt,draft]{amsart}
\usepackage{amscd}
\usepackage{amssymb}
\newtheorem{thm}[equation]{Theorem}
\newtheorem{cor}[equation]{Corollary}
\newtheorem{prop}[equation]{Proposition}
\newtheorem{lem}[equation]{Lemma}
\theoremstyle{definition}
\newtheorem{dfn}[equation]{Definition}
\newtheorem{rem}[equation]{Remark}
\newtheorem{exa}[equation]{Example}

\newtheorem{question}[equation]{Question}
\numberwithin{equation}{section}

\newcommand{\iso}{\stackrel{\simeq}{\rightarrow}}

\newcommand{\surj}{\twoheadrightarrow}
\newcommand{\ar}{\rightarrow}
\newcommand{\opn}{\operatorname}
\newcommand{\cat}[1]{\operatorname{\mathsf{#1}}}
\newcommand{\bdot}{{\textstyle \cdot}}
\newcommand{\ul}{\underline}

\newcommand{\sqbr}[1]{[ \, #1 \, ]}
\newcommand{\mfrak}[1]{\mathfrak{#1}}
\newcommand{\mcal}[1]{\mathcal{#1}}
\newcommand{\msf}[1]{\mathsf{#1}}
\newcommand{\mbf}[1]{\mathbf{#1}}
\newcommand{\mrm}[1]{\mathrm{#1}}
\newcommand{\mbb}[1]{\mathbb{#1}}
\newcommand{\bsym}[1]{\boldsymbol{#1}}
\newcommand{\gfrac}[2]{\genfrac{[}{]}{0pt}{}{#1}{#2}}
\newcommand{\tup}[1]{\textup{#1}}
\renewcommand{\paragraph}[1]{\bigskip \noindent \textbf{#1.}}

\title{Adelic Chern Forms and Applications}
\date{18 March 1998}

\author[Reinhold H\"{u}bl]{Reinhold H\"{u}bl*}
\author[Amnon Yekutieli]{Amnon Yekutieli**}

\address{Fachbereich Mathematik,
Universit\"{a}t Regensburg,
Universit\"{a}tsstrasse 31,
D-93053 Regensburg, GERMANY
\newline \indent
\textit{E-mail address}:
{\tt reinhold.huebl@mathematik.uni-regensburg.de}}

\address{Department of Theoretical Mathematics,
The Weizmann Institute of Science,
Rehovot 76100, ISRAEL
\newline \indent
{\it E-mail address}:
{\tt amnon@wisdom.weizmann.ac.il}}

\thanks{\noindent *Partially supported by the Deutsche
Forschungsgemeinschaft}
\thanks{**Supported by an Allon Fellowship}

\begin{document}

\begin{abstract}
Let $X$ be a variety over a field of characteristic $0$.
Given a vector bundle $E$ on $X$ we construct Chern forms
$c_{i}(E; \nabla) \in \Gamma(X, \mcal{A}^{2i}_{X})$.
Here $\mcal{A}^{\bdot}_{X}$ is the sheaf of {\em Beilinson adeles}
and $\nabla$ is an {\em adelic connection}. When $X$ is smooth
$\mrm{H}^{p} \Gamma(X, \mcal{A}^{\bdot}_{X}) =
\mrm{H}^{p}_{\mrm{DR}}(X)$,
the algebraic De Rham cohomology, and
$c_{i}(E) = [c_{i}(E; \nabla)]$
are the usual Chern classes.

We include three applications of the construction:
(1) existence of adelic secondary (Chern-Simons) characteristic
classes on any smooth $X$ and any vector bundle $E$;
(2) proof of the Bott Residue Formula for a vector field action; and
(3) proof of a Gauss-Bonnet Formula on the level of differential
forms, namely in the De Rham-residue complex.
\end{abstract}

\maketitle

\tableofcontents


\setcounter{section}{-1}
\section{Introduction}

Let $X$ be a scheme of finite type over a field
$k$. According to Beilinson \cite{Be}, given any quasi-coherent
$\mcal{O}_{X}$-module $\mcal{M}$ and an integer $q$, there is a
flasque $\mcal{O}_{X}$-module
$\ul{\mbb{A}}_{\mrm{red}}^{q}(\mcal{M})$,
called the {\em sheaf of adeles}.
This is a generalization of the classical adeles of number theory
(cf.\ Example \ref{exa2.1}).
Moreover, there are homomorphisms
$\partial : \ul{\mbb{A}}_{\mrm{red}}^{q}(\mcal{M}) \ar
\ul{\mbb{A}}_{\mrm{red}}^{q + 1}(\mcal{M})$
which make $\ul{\mbb{A}}_{\mrm{red}}^{\bdot}(\mcal{M})$ into a
complex, and
$\mcal{M} \ar \ul{\mbb{A}}_{\mrm{red}}^{\bdot}(\mcal{M})$
is quasi-isomorphism.

Now let $\Omega_{X / k}^{\bdot}$ be the
algebra of K\"{a}hler differential forms on $X$. In \cite{HY}
we proved that the sheaf
\[ \mcal{A}^{\bdot}_{X} =
\ul{\mbb{A}}_{\mrm{red}}^{\bdot}(\Omega_{X / k}^{\bdot}) =
\bigoplus_{p, q}
\ul{\mbb{A}}_{\mrm{red}}^{q}(\Omega^{p}_{X / k}) \]
is a resolution of $\Omega_{X / k}^{\bdot}$ as differential graded
algebras (DGAs). Therefore when $X$ is smooth,
$\mcal{A}^{\bdot}_{X}$ calculates the
algebraic De Rham cohomology:
$\mrm{H}^{p}_{\mrm{DR}}(X) \cong \mrm{H}^{p} \Gamma(X,
\mcal{A}^{\bdot}_{X})$.
We see that there is an analogy between $\mcal{A}^{\bdot}_{X}$
and the Dolbeault sheaves of smooth forms on a complex-analytic
manifold.

Carrying this analogy further, in this paper we show that when
$\opn{char} k = 0$, any vector bundle $E$ on $X$ admits an
adelic connection $\nabla$. Given such a connection one can
assign adelic Chern forms
$c_{i}(E, \nabla) \in \Gamma(X, \mcal{A}^{2i}_{X})$,
whose classes
$c_{i}(E) := [c_{i}(E, \nabla)] \in \mrm{H}^{2i}_{\mrm{DR}}(X)$
are the usual Chern classes. We include three applications of our
adelic Chern-Weil theory, to demonstrate its effectiveness and
potential.

The idea of using adeles for an algebraic Chern-Weil theory
goes back to Parshin, who constructed a Chern form
$c_{i}(E) \in \ul{\mbb{A}}^{i}(\Omega_{X / k}^{i})$
using an $i$-cocycle on $\mrm{Gl}(\ul{\mbb{A}}^{1}(\mcal{O}_{X}))$
(see \cite{Pa}). Unfortunately we found it quite difficult to
perform calculations with Parshin's forms. Indeed, there is an
inherent complication to any Chern-Weil theory based on
$\mcal{A}^{\bdot}_{X}$.
The DGA $\mcal{A}^{\bdot}_{X}$, with its
Alexander-Whitney product, is not (graded)
commutative. This means that even if one had some kind of
``curvature matrix'' $R$ with entries in $\mcal{A}^{2}_{X}$, one
could not simply evaluate invariant polynomials on $R$.

The problem of noncommutativity was encountered long ago in algebraic
topology, and was dubbed the ``commutative cochain problem''. The
solution, by Thom and Sullivan, was extended to the setup of
cosimplicial DGAs by Bousfield-Gugen\-heim and Hinich-Schechtman (see
\cite{BG}, \cite{HS1}, \cite{HS2}). In our framework this gives a
sheaf of commutative DGAs $\tilde{\mcal{A}}^{\bdot}_{X}$ on $X$,
called the sheaf of {\em Thom-Sullivan adeles}, and a homomorphism of
complexes (``integration on the fiber'')
$\int_{\Delta} :
\tilde{\mcal{A}}^{\bdot}_{X} \ar \mcal{A}^{\bdot}_{X}$.
This map induces an isomorphism of graded algebras
$\mrm{H}^{\bdot}(\int_{\Delta}) :
\mrm{H}^{\bdot} \Gamma(X, \tilde{\mcal{A}}^{\bdot}_{X}) \ar
\mrm{H}^{\bdot} \Gamma(X, \mcal{A}^{\bdot}_{X})$.
We should point out that $\int_{\Delta}$ involves denominators, so
it is necessary to work in characteristic $0$.

Bott discovered a way of gluing together connections defined
locally on a manifold (see \cite{Bo1}). This method was imported to
algebraic geometry by Zhou (in \cite{Zh}), who used \v{C}ech
cohomology. When we tried to write the formulas in terms of adeles, it
became evident that they gave a connection on the Thom-Sullivan adeles
$\tilde{\mcal{A}}^{\bdot}_{X}$. Later we realized that a similar
construction was used by Dupont in the context of simplicial manifolds
(see \cite{Du}).

In the remainder of the Introduction we outline the main results of
our paper.

\paragraph{Adelic Connections}
Let $k$ be a field of characteristic $0$ and $X$ a finite type scheme
over it. The definition of Beilinson adeles on $X$ and their
properties will be reviewed in Section 2.
For now let us just note that the sheaf of adeles
$\tilde{\mcal{A}}^{\bdot}_{X}$ is a commutative DGA, and
$\Omega^{\bdot}_{X / k} \ar \tilde{\mcal{A}}^{\bdot}_{X}$
is a DGA quasi-isomorphism.

Let $\mcal{E}$ be the locally free $\mcal{O}_{X}$-module of rank $r$
associated to the vector bundle $E$.
An {\em adelic connection} on $\mcal{E}$ is by definition a
connection
\[ \nabla : \tilde{\mcal{A}}^{0}_{X} \otimes_{\mcal{O}_{X}}
\mcal{E} \ar
\tilde{\mcal{A}}^{1}_{X} \otimes_{\mcal{O}_{X}} \mcal{E} \]
over the algebra $\tilde{\mcal{A}}^{0}_{X}$.

Such connections are abundant. One way to get an adelic connection is
by choosing, for every point $x$, a basis (or frame; we use these
terms interchangeably)
$\bsym{e}_{x} = (e_{x, 1}, \ldots, e_{x, r})$
for the $\mcal{O}_{X, x}$-module $\mcal{E}_{x}$. We then get a
Levi-Civita connection
\[ \nabla_{x} : \mcal{E}_{x} \ar
\Omega^{1}_{X / k, x} \otimes_{\mcal{O}_{X, x}} \mcal{E}_{x} \]
over the $k$-algebra
$\mcal{O}_{X, x}$.
The Bott gluing mentioned above produces an adelic connection
$\nabla$ (see Proposition \ref{prop3.1}).

\paragraph{Adelic Chern-Weil Homomorphism}
Since $\tilde{\mcal{A}}^{\bdot}_{X}$ is a (graded) commutative DGA,
an adelic connection $\nabla$ on $\mcal{E}$ gives a curvature form
\[ R := \nabla^{2} \in
\Gamma(X, \tilde{\mcal{A}}^{2}_{X} \otimes_{\mcal{O}_{X}}
 \mcal{E}nd(\mcal{E})) . \]

Denote by
$\mrm{S}(\mrm{M}_{r}(k)^{*}) = \mcal{O}(\mrm{M}_{r} \times k)$
the algebra of polynomial functions on $r \times r$ matrices, and let
$I_{r}(k) :=\mrm{S}(\mrm{M}_{r}(k)^{*})^{\mrm{Gl}_{r}(k)}$
be the subalgebra of conjugation-invariant functions.
Denote by $P_{i}$ the $i$-th elementary invariant polynomial, so
$P_{1} = \opn{tr}, \ldots, P_{r} = \opn{det}$.

For any $P \in I_{r}(k)$ the form
$P(R) \in \Gamma(X, \tilde{\mcal{A}}^{\bdot}_{X})$
is closed. So there is a $k$-algebra homomorphism
$w_{\mcal{E}} :
I_{r}(k) \ar
\mrm{H}^{\bdot} \Gamma(X, \tilde{\mcal{A}}^{\bdot}_{X})$,
$P \mapsto [P(R)]$,
called the {\em adelic Chern-Weil homomorphism}.
In Theorem \ref{thm3.2} we prove that $P(\mcal{E}) = w_{\mcal{E}}(P)$
is independent of the connection $\nabla$ (this is true even if $X$
is singular). Defining the $i$-th {\em adelic Chern form} to be
\[ c_{i}(\mcal{E}; \nabla) := \int_{\Delta} P_{i}(R) \in
\Gamma(X, \mcal{A}^{2i}_{X}) , \]
we show the three axioms of Chern classes are satisfied
(Theorem \ref{thm3.4}). Hence when $X$ is smooth over $k$,
\[ c_{i}(\mcal{E}) := [c_{i}(\mcal{E}; \nabla)] \in
\mrm{H}^{2i}_{\mrm{DR}}(X) \]
is the usual $i$-th Chern class.

\paragraph{Secondary Characteristic Classes}
Suppose now that $X$ is a smooth scheme over $k$, and let $\mcal{E}$
be a locally free sheaf of rank $r$ on $X$. Let $P \in I_{r}(k)$ be
an invariant polynomial function of degree $m \geq 2$.
In \cite{BE}, Bloch and Esnault showed that given
an algebraic connection
$\nabla : \mcal{E} \ar \Omega^{1}_{X / k} \otimes_{\mcal{O}_{X}}
\mcal{E}$,
there is a Chern-Simons class
$\mrm{T} P(\mcal{E}, \nabla)$
satisfying
$\mrm{d} \mrm{T} P(\mcal{E}, \nabla) =
P(\mcal{E}) \in \mrm{H}^{2m}_{\mrm{DR}}(X)$.
(We are using the notation of \cite{Es}.)

Since adelic connections always exist, we can construct adelic
secondary characteristic classes on {\em any} smooth $k$-scheme
$X$ and {\em any} locally free sheaf $\mcal{E}$. Theorem \ref{thm4.1}
says that given an adelic connection $\nabla$ there is a class
\[ \mrm{T} P(\mcal{E}; \nabla) \in
\Gamma \left( X, \mcal{A}^{2m-1}_{X} / \mrm{D}(\mcal{A}^{2m-2}_{X})
\right) \]
satisfying
\[ \mrm{D} \mrm{T} P(\mcal{E}; \nabla) = P(\mcal{E}) \in
\mrm{H}^{2m}_{\mrm{DR}}(X) . \]

The existence of adelic secondary characteristic classes, combined
with the action of adeles on the residue complex (see below,
and Theorem \ref{thm6.1}), opens new possibilities for research
on vanishing of cohomology classes (cf.\ \cite{Es}).

\paragraph{Bott Residue Formula}
The adeles of differential forms can be integrated.
If $\opn{dim} X = n$, each maximal chain of points
$\xi = (x_{0}, \ldots, x_{n})$ determines a local integral
$\opn{Res}_{\xi} : \Gamma(X, \mcal{A}^{2n}_{X}) \ar k$
(cf.\ \cite{Be}, \cite{Ye1}). If $X$ is smooth and proper then
the global map
\[ \int_{X}  := \sum_{\xi} \opn{Res}_{\xi} :
\mrm{H}^{2n}_{\mrm{DR}}(X) =
\mrm{H}^{2n} \Gamma(X, \mcal{A}^{\bdot}_{X}) \ar k \]
coincides with the usual ``algebraic integral'' of, say, \cite{Ha1}.

Assume $X$ is a smooth projective variety of dimension $n$.
Let $P \in I_{r}(k)$ be a homogeneous polynomial of degree $n$,
so that $P = Q(P_{1}, \ldots, P_{r})$ for some polynomial $Q$ in $r$
variables. Let $v \in \Gamma(X, \mcal{T}_{X})$ be a vector field
with isolated zeroes, and assume $v$ acts on the locally free sheaf
$\mcal{E}$. For each zero $z$ of $v$ there is a local
invariant $P(v, \mcal{E}, z) \in k$, which has an explicit expression
in terms of local coordinates. Theorem \ref{thm5.1} says that
\[ \int_{X} Q(c_{1}(\mcal{E}), \ldots, c_{r}(\mcal{E})) =
\sum_{v(z) = 0} P(v, \mcal{E}, z) . \]
The proof of the theorem follows the steps of Bott's original
proof in \cite{Bo2}, translated to adeles and algebraic residues.
Example \ref{exa5.1} provides an explicit illustration of the result
in the case of a nonreduced zero $z$.

We should of course mention the earlier algebraic proof of the Bott
Residue Formula for isolated zeroes, by Carrell-Lieberman \cite{CL},
which uses Grothendieck's global duality.

There is also a Bott Residue Formula for group actions, which is best
stated as a localization formula in equivariant cohomology (cf.\
\cite{AB}). Recently this formula was used in enumerative geometry,
see for instance \cite{ES} and \cite{Ko}.
Edidin-Graham \cite{EG} proved Bott's formula in the equivariant
intersection ring.

\paragraph{The Gauss-Bonnet Formula}
Let $k$ be a perfect field of any characteristic, and let $X$ be a
finite type $k$-scheme. The residue complex $\mcal{K}^{\bdot}_{X}$ is
by definition the Cousin complex of $\pi^{!} k$, where
$\pi : X \ar \opn{Spec} k$ is the structural morphism (cf.\
\cite{RD}). Each $\mcal{K}^{q}_{X}$ is a quasi-coherent sheaf. Let
\[ \mcal{F}^{\bdot}_{X} := \mcal{H}om_{\mcal{O}_{X}}(
\Omega^{\bdot}_{X / k}, \mcal{K}^{\bdot}_{X}) \]
which is a graded sheaf in the obvious way.
According to \cite{EZ} or \cite{Ye3}, there is an operator
$\mrm{D} : \mcal{F}^{i}_{X} \ar \mcal{F}^{i + 1}_{X}$
which makes $\mcal{F}^{\bdot}_{X}$ into a DG
$\Omega^{\bdot}_{X / k}$-module. $\mcal{F}^{\bdot}_{X}$ is called the
{\em De Rham-residue complex}.
When $X$ is smooth,
$\mrm{H}_{i}^{\mrm{DR}}(X) =
\mrm{H}^{-i} \Gamma(X, \mcal{F}^{\bdot}_{X})$.

In \cite{Ye5} it is we proved that there is a natural
structure of right DG
$\mcal{A}^{\bdot}_{X}$-module on $\mcal{F}^{\bdot}_{X}$, extending
the $\Omega^{\bdot}_{X / k}$-module structure
(cf.\ Theorem \ref{thm6.1}).
The action is ``by taking residues.''
When $f : X \ar Y$ is proper then
$\opn{Tr}_{f} : f_{*} \mcal{F}^{\bdot}_{X} \ar \mcal{F}^{\bdot}_{Y}$
is a homomorphism of DG $\mcal{A}^{\bdot}_{Y}$-modules.
If we view the adeles $\mcal{A}^{p, q}_{X}$ as an algebraic analog of
the smooth forms of type $(p, q)$ on a complex-analytic manifold,
then $\mcal{F}^{-p, -q}_{X}$ is the analog of the currents of
type $(p, q)$.

Suppose $\opn{char} k = 0$, $X$ is smooth irreducible of dimension
$n$, $\mcal{E}$ is a locally free $\mcal{O}_{X}$-module of rank $r$,
$v$ is a regular section of $\mcal{E}$ and $Z$ is its zero scheme.
Let $C_{X}, C_{Z} \in \Gamma(X, \mcal{F}^{\bdot}_{X})$ be the
fundamental classes. In Theorem \ref{thm7.1} we prove the following
version of the Gauss-Bonnet Formula:
there is an adelic connection $\nabla$ on $\mcal{E}$ satisfying
\[ C_{X} \cdot c_{r}(\mcal{E}, \nabla) =
(-1)^{m} C_{Z} \in \mcal{F}_{X}^{-2(n - r)}  \]
with $m = nr + \binom{r+1}{2}$.

Observe that this formula is on the level of differential forms.
Passing to (co)ho\-mo\-logy we recover the familiar formula
$c_{r}(\mcal{E}) \smile [X] = [Z] \in
\mrm{H}_{2n - r}^{\mrm{DR}}(X)$
(cf.\ \cite{Fu} \S 14.1 ).

\paragraph{Acknowledgements}
The authors wish to thank S.\ Kleiman who suggested studying the
adelic approach to the Bott Residue Formula, J.\ Lipman for calling
our attention to Zhou's work, V.\ Hinich for explaining to us the
construction of the Thom-Sullivan cochains, D.\ Blanc for help with
cosimplicial groups, and P.\ Sastry for discussions. Special thanks
also go to P.\ Golginger and W.\ Krinner for collecting the available
material on algebraic connections (\cite{Go}, \cite{Kr}). Part of the
research was done while the first author was visiting the Weizmann
Institute. He wishes to express his gratitude towards this institution
for its hospitality during his stay there.

\section{Cosimplicial Algebras and their Normalizations}

In this section we review some well known facts about cosimplicial
objects, and also discuss the less known Thom-Sullivan normalization.
Our sources are \cite{Ma}, \cite{ML}, \cite{BG} and \cite{HS1}.

Denote by $\Delta$ the category whose objects are
the finite ordered sets $[n] := \{ 0, \ldots, n \}$, and whose
morphisms are the monotone nondecreasing functions $[m] \ar [n]$. Let
$\partial^{i} : [n-1] \ar [n]$ stand for the $i$-th coface map, and
let $s^{i} : [n+1] \ar [n]$ stand for the $i$-th codegeneracy map.
These maps generate by composition all morphisms in $\Delta$.
Let $\Delta^{\circ}$ denote the opposite category.
By definition a {\em simplicial object} in a category $\msf{C}$ is a
functor $S : \Delta^{\circ} \ar \msf{C}$. Often one writes $S_{n}$
instead of the object $S[n] \in \msf{C}$. A {\em cosimplicial object}
is a functor $S : \Delta \ar \msf{C}$.
Denote by $\Delta^{\circ} \msf{C}$ (resp.\ $\Delta \msf{C}$)
the category of simplicial (resp.\ cosimplicial) objects in $\msf{C}$.

\begin{exa}
Let $P$ be a partially ordered set.
A simplex (or chain) of length $n$ in $P$ is a sequence
$\sigma = (x_{0}, \ldots, x_{n})$, with $x_{i} \leq x_{i+1}$.
More generally, if $P$ is a category, then an $n$-simplex is a functor
$\sigma : [n] \ar P$. Letting
$S(P)_{n}$ be the set of $n$-simplices in $P$, we see that
$S(P)$ is a simplicial set.
\end{exa}

\begin{exa} \label{exa1.5}
If we take $P = [n]$, then we get the standard simplicial complex
$\Delta^{n} \in \Delta^{\circ} \msf{Sets}$.
As a functor $\Delta^{\circ} \ar \msf{Sets}$ one has
$\Delta^{n} = \opn{Hom}_{\Delta}(-, [n])$.
Observe that
\[ \Delta_{m}^{n} = \opn{Hom}_{\Delta}([m], [n]) =
\{ (i_{0}, \ldots, i_{m})\ |\
0 \leq i_{0} \leq \cdots \leq i_{m} \leq n
\} . \]
\end{exa}

\begin{exa} \label{exa1.1}
Given a scheme $X$, specialization defines a partial ordering on its
underlying set of points: $x \leq y$ if $y \in \overline{\{ x \}}$.
We denote by $S(X)$ the resulting simplicial set.
\end{exa}

\begin{exa} \label{exa1.2}
Let $\Delta^{n}_{\mrm{top}}$ be the standard realization of
$\Delta^{n}$, i.e.\ the compact topological space
\[ \{ (t_{0}, \ldots, t_{n})\ |\ t_{i} \geq 0 \text{ and }
\sum t_{i} =1 \} \subset \mbb{R}^{n+1} . \]
Then $\Delta_{\mrm{top}} = \{ \Delta^{n}_{\mrm{top}} \}$
is a cosimplicial topological space.
\end{exa}

Let $k$ be any commutative ring. By a differential graded (DG)
$k$-module we mean a cochain complex, namely a graded module
$M = \bigoplus_{q \in \mbb{Z}} M^{q}$
with an endomorphism $\mrm{d}$ of degree $1$ satisfying
$\mrm{d}^{2} = 0$.
By a differential graded algebra (DGA) over $k$ we mean a DG module
$A$ with a DG homomorphism $A \otimes_{k} A \ar A$. So $A$ is
neither assumed to be commutative nor associative.

\begin{exa} \label{exa1.3}
Let $t_{0}, t_{1}, \ldots$ be indeterminates (of degree $0$). Define
\[ R_{n} := k \sqbr{ t_{0}, \ldots, t_{n} } /
(t_{0} + \cdots + t_{n} -1) \]
and
$\Delta^{n}_{k} := \opn{Spec} R_{n}$.
Then  as in the previous example, $\Delta_{k} = \{ \Delta^{n}_{k} \}$
is a cosimplicial scheme. Letting
$\Omega^{\bdot}(\Delta^{n}_{k}) := \Omega^{\bdot}_{R_{n} / k}$,
we see that
$\Omega^{\bdot}(\Delta_{k}) := \{
\Omega^{\bdot}(\Delta^{n}_{k}) \}$
is a simplicial DGA over $k$.
\end{exa}

Consider a cosimplicial $k$-module
$M = \{ M^{q} \} \in \Delta \msf{Mod}(k)$.
Its standard normalization is the DG module
$(\mrm{N} M, \partial)$ whose degree $q$ piece is
$\mrm{N}^{q} M := \bigcap \opn{Ker}(s^{i}) \subset M^{q}$,
and $\partial := \sum (-1)^{i} \partial^{i}$.
Now suppose
$M = \{ M^{\bdot, q} \} \in \Delta \msf{DGMod}(k)$,
i.e.\ a cosimplicial DG $k$-module. Each
$M^{\bdot, q} = \bigoplus_{p \in \mbb{Z}} M^{p, q}$
is a DG module with operator
$\mrm{d} : M^{p, q} \ar M^{p+1, q}$, and each
$M^{p,\bdot}$ is a cosimplicial $k$-module.
Define $\mrm{N}^{p,q} M := \mrm{N}^{q} M^{p,\bdot}$,
$\mrm{N}^{i} M := \bigoplus_{p+q=i} \mrm{N}^{p,q} M$,
and
$\mrm{N} M := \bigoplus_{i} \mrm{N}^{i} M$.
Then  $\mrm{N} M$ is a DG module with coboundary operator
$\mrm{D} := \mrm{D}' + \mrm{D}''$,
where
$\mrm{D}' := (-1)^{q} \mrm{d} : \mrm{N}^{p,q} M \ar
\mrm{N}^{p+1,q} M$
and
$\mrm{D}'' := \partial : \mrm{N}^{p,q} M \ar \mrm{N}^{p,q+1} M$.
Another way to visualize this is by defining for each $q$ a DG module
$\mrm{N}^{\bdot,q} M := \bigcap \opn{Ker}(s^{i})
\subset M^{\bdot,q}[-q]$
(the shift by $-q$), so the operator is indeed $\mrm{D}'$.
Then
$\mrm{D}'' = \partial : \mrm{N}^{\bdot,q} M \ar
\mrm{N}^{\bdot,q+1} M$
has degree $1$ and
$\mrm{N} M = \bigoplus_{q} \mrm{N}^{\bdot,q} M$.

If $A$ is a cosimplicial DGA, that is $A \in \Delta \cat{DGA}(k)$,
then $\mrm{N} A$ is a DGA with
the Alexander-Whitney product. For any
$a \in \mrm{N}^{p, q} A$ and $b \in \mrm{N}^{p',q'} A$ one has
\begin{equation} \label{eqn1.2}
a \cdot b =
\partial^{-}(a) \cdot \partial^{+}(b) \in
\mrm{N}^{p+p', q+q'} A ,
\end{equation}
where
$\partial^{-} : [q] \ar [q + q']$ is the simplex $(0,1, \ldots, q)$
and
$\partial^{+} : [q'] \ar [q + q']$ is the simplex
$(q,q+1, \ldots, q+q')$ (cf.\ Example \ref{exa1.5}).
Note that if each algebra $A^{\bdot, q}$ is associative, then so is
$\mrm{N} A$; however $\mrm{N} A$ is usually not commutative.
If $M$ is a cosimplicial DG left $A$-module then $\mrm{N} M$ is a DG
left $\mrm{N} A$-module.

We shall need another normalization of cosimplicial objects.
The definition below is extracted from the work of Hinich-Schechtman,
cf.\ \cite{HS1}, \cite{HS2}. Fix a commutative
$\mbb{Q}$-algebra $k$.

\begin{dfn} \label{dfn1.1}
Suppose $M= \{ M^{q} \}$ is a cosimplicial $k$-module.
Let
\[ \tilde{\mrm{N}}^{q} M \subset
\prod_{l=0}^{\infty} \left( \Omega^{q}(\Delta^{l}_{\mbb{Q}})
\otimes_{\mbb{Q}} M^{l} \right) \]
be the submodule consisting of all elements
$\bsym{u} = (u_{0}, u_{1}, \ldots)$,
$u_{l} \in \Omega^{q}(\Delta^{l}_{\mbb{Q}}) \otimes_{\mbb{Q}} M^{l}$,
s.t.\
\begin{eqnarray} \label{eqn1.5}
(1 \otimes \partial^{i}) u_{l} & = & (\partial_{i} \otimes 1)
u_{l+1} \\
(s_{i} \otimes 1) u_{l} & = & (1 \otimes s^{i}) u_{l+1}
\end{eqnarray}
for all $0 \leq l$, $0 \leq i \leq l+1$.
Given a cosimplicial DG $k$-module $M = \{ M^{\bdot, q} \}$, let
$\tilde{\mrm{N}}^{p,q} M := \tilde{\mrm{N}}^{q} M^{p,\bdot}$,
$\tilde{\mrm{N}}^{i} M := \bigoplus_{p+q=i}
\tilde{\mrm{N}}^{p,q} M$
and
$\tilde{\mrm{N}} M := \bigoplus_{i} \tilde{\mrm{N}}^{i} M$.
Define
$\mrm{D}' := (-1)^{q} \otimes \mrm{d}$,
$\mrm{D}'' := \mrm{d} \otimes 1$ and
$\mrm{D} := \mrm{D}' + \mrm{D}''$.
The resulting complex $(\tilde{\mrm{N}} M, \mrm{D})$
is called the {\em Thom-Sullivan normalization} of $M$.
If $A$ is a cosimplicial $k$-DGA, then $\tilde{\mrm{N}} A$
inherits the component-wise multiplication from
$\prod_{l} \Omega^{q}(\Delta^{l}_{\mbb{Q}})
\otimes_{\mbb{Q}} A^{p, l}$,
so it is a DGA.
\end{dfn}

In the definition above the signs are in agreement with the usual
conventions; keep in mind that $M^{p,q}$ is in degree $p$
(cf.\ \cite{ML} Ch.\ VI \S 7). It is clear that if each DGA
$A^{\bdot, q}$
is commutative (resp.\ associative), then so is $\tilde{\mrm{N}} A$.

Usual integration on the real simplex $\Delta^{l}_{\mrm{top}}$
yields a $\mbb{Q}$-linear map of degree $0$,
$\int_{\Delta^{l}} : \Omega^{\bdot}(\Delta^{l}_{\mbb{Q}})
\ar \mbb{Q}[-l]$,
such that
$\int_{\Delta^{l}}(\mrm{d} t_{1} \wedge \cdots \wedge
\mrm{d} t_{l}) = \frac{1}{l !}$.
By linearity, for any cosimplicial DG module $M$ this extends to
a degree $0$ homomorphism
\[ \int_{\Delta^{l}} : \Omega^{\bdot}(\Delta^{l}_{\mbb{Q}})
\otimes_{\mbb{Q}} M^{\bdot,l} \ar
\mbb{Q}[-l] \otimes_{\mbb{Q}} M^{\bdot,l} = M^{\bdot,l}[-l] . \]
Note that $\int_{\Delta^{l}}$ sends
$\mrm{D}' := (-1)^{q} \otimes \mrm{d}$ to
$\mrm{D}' := (-1)^{l} \mrm{d}$. Define
$\int_{\Delta} : \tilde{\mrm{N}} M \ar \bigoplus_{q} M^{\bdot,q}[-q]$
by:
\begin{equation} \label{eqn1.3}
\int_{\Delta} (u_{0}, u_{1}, \ldots) :=
\sum_{l = 0}^{\infty} \int_{\Delta^{l}} u_{l} .
\end{equation}

\begin{lem}
The image of $\int_{\Delta}$ lies inside $\mrm{N} M$, and
$\int_{\Delta} : \tilde{\mrm{N}} M \ar \mrm{N} M$
is a $k$-linear homomorphism of complexes.
\end{lem}

\begin{proof} A direct verification, amounting to Stoke's Theorem on
$\Delta^{l}_{\mrm{top}}$.
\end{proof}

What we have is a natural transformation
$\int_{\Delta} : \tilde{\mrm{N}} \ar \mrm{N}$
of functors
$\Delta \msf{DGMod}(k) \ar$ \linebreak
$\msf{DGMod}(k)$.

\begin{thm} \label{thm1.1} \textup{(Simplicial De Rham Theorem)}\
Let $M$ a cosimplicial DG $\mbb{Q}$-module. Then
$\int_{\Delta} : \tilde{\mrm{N}} M \ar \mrm{N} M$
is a quasi-isomorphism.
\end{thm}

\begin{thm} \label{thm1.2}
Let $k$ be a commutative $\mbb{Q}$-algebra and $A$ a cosimplicial DG
$k$-algebra. Then
$\mrm{H}(\int_{\Delta}) : \mrm{H} \tilde{\mrm{N}} A \ar
\mrm{H} \mrm{N} A$
is an isomorphism of graded $k$-algebras. If $M$ is a cosimplicial
DG $A$-module, then
$\mrm{H}(\int_{\Delta}) : \mrm{H} \tilde{\mrm{N}} M \ar
\mrm{H} \mrm{N} M$
is an isomorphism of graded $\mrm{H} \tilde{\mrm{N}} A$-modules.
\end{thm}

The proof are essentially contained in \cite{BG}  and \cite{HS1}.
For the sake of completeness we include proofs in Appendix A
of this paper.

\section{Adeles of Differential Forms}

In this section we apply the constructions of Section 1 to the
cosimplicial DGA $\ul{\mbb{A}}(\Omega^{\bdot}_{X/k})$ on a scheme $X$.
This will give two DGAs, $\mcal{A}^{\bdot}_{X}$ and
$\tilde{\mcal{A}}^{\bdot}_{X}$, which are resolutions of
$\Omega^{\bdot}_{X/k}$.

Let us begin with a review of {\em Beilinson adeles} on a noetherian
scheme $X$ of finite dimension. A chain of points in $X$ is a sequence
$\xi = (x_{0}, \ldots, x_{q})$ of points with
$x_{i+1} \in \overline{ \{ x_{i} \} }$. Denote by $S(X)_{q}$ the set
of length $q$ chains, so $\{ S(X)_{q} \}_{q \geq 0}$ is a simplicial
set.
For $T \subset S(X)_{q}$ and $x \in X$ let
\[ \hat{x} T := \{ (x_{1}, \ldots, x_{q}) \mid
(x, x_{1}, \ldots, x_{q}) \in T \} . \]
According to \cite{Be} there is a unique collection of functors
$\mbb{A}(T, -) : \msf{Qco}(X) \ar \msf{Ab}$,
indexed by $T \subset S(X)_{q}$, each of which commuting with direct
limits, and satisfying
\[ \mbb{A}(T, \mcal{M}) = \begin{cases}
\prod_{x \in X} \lim_{\leftarrow n} \mcal{M}_{x} / \mfrak{m}_{x}^{n}
\mcal{M}_{x} & \text{if } q = 0 \\[1mm]
\prod_{x \in X} \lim_{\leftarrow n}
\mbb{A}(\hat{x} T, \mcal{M}_{x} / \mfrak{m}_{x}^{n} \mcal{M}_{x})
& \text{if } q > 0
\end{cases} \]
for $\mcal{M}$ coherent.
Here $\mfrak{m} \subset \mcal{O}_{X, x}$ is the maximal ideal
and $\mcal{M}_{x} / \mfrak{m}_{x}^{n} \mcal{M}_{x}$ is treated as a
quasi-coherent sheaf with support $\overline{ \{ x \} }$.
Furthermore each $\mbb{A}(T, -)$ is exact.

For a single chain $\xi$ one also writes
$\mcal{M}_{\xi} := \mbb{A}(\{ \xi \}, \mcal{M})$,
and this is the {\em Beilinson completion} of $\mcal{M}$ along $\xi$.
Then
\begin{equation} \label{eqn1.6}
\mbb{A}(T, \mcal{M}) \subset \prod_{\xi \in T} \mcal{M}_{\xi}
\end{equation}
which permits us to consider the adeles as a ``restricted product.''
For $q = 0$ and $\mcal{M}$ coherent we have
$\mcal{M}_{(x)} = \widehat{\mcal{M}}_{x}$, the $\mfrak{m}_{x}$-adic
completion, and (\ref{eqn1.6}) is an equality.
In view of this we shall say that $\mbb{A}(T, \mcal{M})$ is the
group of adeles combinatorially supported on $T$ and with values in
$\mcal{M}$.

Define a presheaf
$\ul{\mbb{A}}(T, \mcal{M})$ by
\begin{equation}
\Gamma(U, \ul{\mbb{A}}(T, \mcal{M})) :=
\mbb{A}(T \cap S(U)_{q}, \mcal{M})
\end{equation}
for $U \subset X$ open. Then $\ul{\mbb{A}}(T, \mcal{M})$ is a flasque
sheaf. Also $\ul{\mbb{A}}(T, \mcal{O}_{X})$ is a flat
$\mcal{O}_{X}$-algebra, and
$\ul{\mbb{A}}(T, \mcal{M}) \cong \ul{\mbb{A}}(T, \mcal{O}_{X})
\otimes_{\mcal{O}_{X}} \mcal{M}$.

For every $q$ define the sheaf of degree $q$ {\em Beilinson adeles}
\[ \ul{\mbb{A}}^{q}(\mcal{M}) :=
\ul{\mbb{A}}(S(X)_{q}, \mcal{M}) . \]
Then
$\ul{\mbb{A}}(\mcal{M}) =
\{ \ul{\mbb{A}}^{q}(\mcal{M}) \}_{q \in \mbb{N}}$
is a cosimplicial sheaf.
The standard normalization
$\mrm{N}^{q} \ul{\mbb{A}}(\mcal{M})$
is canonically isomorphic to the sheaf
$\ul{\mbb{A}}^{q}_{\mrm{red}}(\mcal{M}) :=
\ul{\mbb{A}}(S(X)_{q}^{\mrm{red}}, \mcal{M})$,
where $S(X)_{q}^{\mrm{red}}$ is the set of nondegenerate chains.
Note that
$\ul{\mbb{A}}^{q}_{\mrm{red}}(\mcal{M}) = 0$
for all $q > \opn{dim} X$.
A fundamental theorem of Beilinson says that the canonical
homomorphism
$\mcal{M} \ar \ul{\mbb{A}}^{\bdot}_{\mrm{red}}(\mcal{M})$
is a quasi-isomorphism. We see that
$\mrm{H}^{q} \Gamma(X, \ul{\mbb{A}}^{\bdot}_{\mrm{red}}(\mcal{M}))
= \mrm{H}^{q}(X, \mcal{M})$.

The complex
$\ul{\mbb{A}}^{\bdot}_{\mrm{red}}(\mcal{O}_{X})$
is a DGA, with the Alexander-Whitney product. For local sections
$a \in \ul{\mbb{A}}^{q}_{\mrm{red}}(\mcal{O}_{X})$
and
$b \in \ul{\mbb{A}}^{q'}_{\mrm{red}}(\mcal{O}_{X})$
the product is
$a \cdot b = \partial^{-}(a) \cdot \partial^{+}(b) \in
\ul{\mbb{A}}^{q+q'}_{\mrm{red}}(\mcal{O}_{X})$,
where $\partial^{-}$ and $\partial^{+}$ correspond respectively to the
initial and final segments of $(0, \ldots, q, \ldots, q+q')$.
This algebra is not (graded) commutative.
For proofs and more details turn to \cite{Hr}, \cite{Ye1} Chapter 3
and \cite{HY} Section 1.

\begin{exa} \label{exa2.1}
Suppose $X$ is a nonsingular curve. The relation to the classical
ring of adeles $\mbb{A}(X)$ of Chevalley and Weil is
$\mbb{A}(X) =
\Gamma(X, \ul{\mbb{A}}^{1}_{\mathrm{red}}(\mcal{O}_{X}))$.
\end{exa}

Now assume $X$ is a finite type scheme over the noetherian ring $k$.
In \cite{HY} it was shown that given any differential operator
(DO) $D : \mcal{M} \ar \mcal{N}$ there is an induced DO
$D : \ul{\mbb{A}}^{q}(\mcal{M}) \ar
\ul{\mbb{A}}^{q}(\mcal{N})$.
Applying this to the De Rham complex $\Omega^{\bdot}_{X/k}$
we get a cosimplicial DGA
$\ul{\mbb{A}}(\Omega^{\bdot}_{X/k})$. The {\em De Rham adele complex}
is the DGA
\[ \mcal{A}^{\bdot}_{X} := \mrm{N}
\ul{\mbb{A}}(\Omega^{\bdot}_{X/k}) . \]
Since
\[ \mcal{A}^{p, q}_{X}
\cong \ul{\mbb{A}}^{q}_{\mrm{red}}(\mcal{O}_{X})
\otimes_{\mcal{O}_{X}} \Omega^{p}_{X/k}  \]
we see that $\mcal{A}^{\bdot}_{X}$ is bounded.
By a standard double complex spectral sequence argument (see
\cite{HY} Proposition 2.1) we get:

\begin{prop} \label{prop2.1}
The natural DGA map
$\Omega_{X / k}^{\bdot} \ar \mcal{A}_{X}^{\bdot}$
is a quasi-iso\-morph\-ism of sheaves. Hence
$\mrm{H}^{\bdot}(X, \Omega_{X / k}^{\bdot}) \cong
\mrm{H}^{\bdot} \Gamma(X, \mcal{A}_{X}^{\bdot})$.
\end{prop}

Let us examine the DGA $\mcal{A}^{\bdot}_{X}$ a little closer.
The operators are
$\mrm{D}' = (-1)^{q} \mrm{d} : \mcal{A}^{p, q}_{X} \ar
\mcal{A}^{p+1, q}_{X}$,
$\mrm{D}'' =  \partial : \mcal{A}^{p, q}_{X} \ar
\mcal{A}^{p, q+1}_{X}$
and
$\mrm{D} = \mrm{D}' + \mrm{D}''$.
As for the multiplication, consider local sections
$a  \in \ul{\mbb{A}}_{\mrm{red}}^{q}(\mcal{O}_{X})$,
$b  \in \ul{\mbb{A}}_{\mrm{red}}^{q'}(\mcal{O}_{X})$,
$\alpha \in \Omega^{p}_{X/k}$ and
$\beta \in \Omega^{p'}_{X/k}$. Then
\begin{equation} \label{eqn2.3}
(a \otimes \alpha) \cdot (b \otimes \beta) =
(-1)^{q'p} \partial^{-}(a) \cdot \partial^{+}(b)
\otimes  \alpha \wedge \beta
\in \mcal{A}_{X}^{p+p',q+q'}
\end{equation}
(cf.\ formula (\ref{eqn1.2})).

\begin{rem}
In the analogy to sheaves of smooth forms on a complex-analytic
manifold, our operators $\mrm{D}'$, $\mrm{D}''$ play the roles of
$\partial$, $\bar{\partial}$ respectively. Note however that here
$\mcal{A}_{X}^{p,q}$ is not a locally free $\mcal{A}_{X}^{0}$-module
for $0 < q \leq n$, even when $X$ is smooth. The
same is true also for the sheaves
$\tilde{\mcal{A}}_{X}^{p, q}$ defined below.
\end{rem}

If $k$ is a perfect field and $X$ is an integral scheme of dimension
$n$, then each maximal chain $\xi = (x_{0}, \ldots, x_{n})$ defines a
$k$-linear map
$\opn{Res}_{\xi} : \Omega^{n}_{X / k} \ar k$ called the {\em Parshin
residue} (cf.\ \cite{Ye1} Definition 4.1.3). By (\ref{eqn1.6}) we
obtain $\opn{Res}_{\xi} : \Gamma(X, \mcal{A}^{2n}_{X}) \ar k$.

\begin{prop} \label{prop2.4}
Suppose $k$ is a perfect field.
\begin{enumerate}
\item Given $\alpha \in \Gamma(X, \mcal{A}^{2n}_{X})$, one has
$\opn{Res}_{\xi} \alpha = 0$ for all but finitely many
$\xi$. Hence
$\int_{X} := \sum_{\xi} \opn{Res}_{\xi} :
\Gamma(X, \mcal{A}^{2n}_{X}) \ar k$
is well-defined.
\item If $X$ is proper then $\int_{X} \mrm{D} \beta = 0$ for all
$\beta \in \Gamma(X, \mcal{A}^{2n - 1}_{X})$. Hence
$\int_{X} : \mrm{H}^{2n}(X, \Omega^{\bdot}_{X / k}) \ar k$
is well-defined.
\item If $X$ is smooth and proper then
$\int_{X} : \mrm{H}^{2n}_{\mrm{DR}}(X) \ar k$
coincides with the nondegenerate map of \cite{Ha1}.
\end{enumerate}
\end{prop}

\begin{proof}
1.\ See \cite{HY} Proposition 3.4.\\
2.\ This follows from the Parshin-Lomadze Residue Theorem (\cite{Ye1}
Theorem 4.2.15).\\
3.\ By \cite{HY} Theorem 3.1 and \cite{Ye3} Corollary 3.8.
\end{proof}

Now assume $k$ is any $\mbb{Q}$-algebra.
The Thom-Sullivan normalization determines a sheaf
$\tilde{\mrm{N}}^{q} \ul{\mbb{A}}(\mcal{M})$, where
$\Gamma(U, \tilde{\mrm{N}}^{q} \ul{\mbb{A}}(\mcal{M})) =
\tilde{\mrm{N}}^{q} \Gamma(U, \ul{\mbb{A}}(\mcal{M}))$.
Applying this to the cosimplicial DGA
$\ul{\mbb{A}}(\Omega^{\bdot}_{X/k})$
we obtain:

\begin{dfn}
The sheaf of {\em Thom-Sullivan adeles} is the sheaf of DGAs
\[ \tilde{\mcal{A}}^{\bdot}_{X} := \tilde{\mrm{N}}
\ul{\mbb{A}}(\Omega^{\bdot}_{X/k}) . \]
\end{dfn}

$(\tilde{\mcal{A}}^{\bdot}_{X}, \mrm{D})$ is an associative,
commutative DGA. The natural map
$\Omega^{\bdot}_{X/k} \ar \tilde{\mcal{A}}^{\bdot}_{X}$
is an injective DGA homomorphism.
The coboundary operator on $\tilde{\mcal{A}}^{\bdot}_{X}$ is
$\mrm{D} = \mrm{D}' + \mrm{D}''$,
where
$\mrm{D}' : \tilde{\mcal{A}}^{p,q}_{X} \ar
\tilde{\mcal{A}}^{p+1, q}_{X}$
and
$\mrm{D}'' : \tilde{\mcal{A}}^{p,q}_{X} \ar
\tilde{\mcal{A}}^{p, q+1}_{X}$.
The ``integral on the fibers'' $\int_{\Delta}$
sheafifies to give a degree $0$ DG $\Omega^{\bdot}_{X/k}$-module
homomorphism
$\int_{\Delta} : \tilde{\mcal{A}}^{\bdot}_{X} \ar
\mcal{A}^{\bdot}_{X}$.
This is not an algebra homomorphism! However:

\begin{prop} \label{prop2.2}
For every open set $U \subset X$,
$\mrm{H}^{\bdot}(\int_{\Delta}) :
\mrm{H}^{\bdot} \Gamma(U, \tilde{\mcal{A}}^{\bdot}_{X}) \ar
\mrm{H}^{\bdot} \Gamma(U, \mcal{A}^{\bdot}_{X})$
is an isomorphism of graded $k$-algebras.
\end{prop}

\begin{proof}
Apply Theorems \ref{thm1.1} and \ref{thm1.2} to the cosimplicial DGA
$\Gamma(U, \ul{\mbb{A}}(\Omega^{\bdot}_{X/k}))$.
\end{proof}

\begin{rem}
We do not know whether the sheaves $\tilde{\mcal{A}}^{p, q}_{X}$
are flasque. The DGA $\tilde{\mcal{A}}^{\bdot}_{X}$ is not bounded;
however, letting
$n := \sup \{ p \mid \Omega^{p}_{X / k, x} \neq 0 \text{ for some }
x \in X \}$, then $n < \infty$ and
$\tilde{\mcal{A}}^{p, q}_{X} = 0$ for all $p > n$.
\end{rem}

\begin{cor} \label{cor2.2}
If $X$ is smooth over $k$, then the homomorphisms
$\Omega^{\bdot}_{X/k} \ar \tilde{\mcal{A}}^{\bdot}_{X} \ar
\mcal{A}^{\bdot}_{X}$
induce isomorphisms of graded $k$-algebras
\[ \mrm{H}^{\bdot}_{\mrm{DR}}(X) =
\mrm{H}^{\bdot}(X, \Omega^{\bdot}_{X / k}) \cong
\mrm{H}^{\bdot} \Gamma(X, \tilde{\mcal{A}}^{\bdot}_{X}) \cong
\mrm{H}^{\bdot} \Gamma(X, \mcal{A}^{\bdot}_{X}) . \]
\end{cor}

Given a quasi-coherent sheaf $\mcal{M}$ set
\begin{equation} \label{eqn2.4}
\tilde{\mcal{A}}^{p,q}_{X}(\mcal{M}) :=
\tilde{\mrm{N}}^{q} \ul{\mbb{A}}(\Omega^{p}_{X/k}
\otimes_{\mcal{O}_{X}} \mcal{M}) .
\end{equation}
In particular we have
$\tilde{\mcal{A}}^{p,q}_{X} =
\tilde{\mcal{A}}^{0,q}_{X}(\Omega^{p}_{X/k})$.

\begin{lem} \label{lem2.1}
\begin{enumerate}
\item Let $\mcal{M}$ be a quasi-coherent sheaf. Then the complex
\[ 0 \ar \mcal{M} \ar \tilde{\mcal{A}}^{0,0}_{X}(\mcal{M})
\xrightarrow{\mrm{D}''} \tilde{\mcal{A}}^{0,1}_{X}(\mcal{M})
\xrightarrow{\mrm{D}''} \cdots \]
is exact.
\item If $\mcal{E}$ is locally free of finite rank, then
$\tilde{\mcal{A}}^{p,q}_{X}(\mcal{E}) \cong
\tilde{\mcal{A}}^{p,q}_{X} \otimes_{\mcal{O}_{X}} \mcal{E}$.
\item Suppose $\mrm{d} : \mcal{M} \ar \mcal{N}$ is a $k$-linear DO.
Then $\mrm{d}$ extends to a DO
$\mrm{d} : \tilde{\mcal{A}}^{0,q}_{X}(\mcal{M}) \ar
\tilde{\mcal{A}}^{0,q}_{X}(\mcal{N})$
which commutes with $\mrm{D}''$.
\end{enumerate}
\end{lem}

\begin{proof}
1.\ Use the quasi-isomorphism
\[ \int_{\Delta} : \tilde{\mcal{A}}^{0,\bdot}_{X}(\mcal{M}) =
\tilde{\mrm{N}} \ul{\mbb{A}}(\mcal{M}) \ar
\mrm{N} \ul{\mbb{A}}(\mcal{M}) =
\ul{\mbb{A}}^{\bdot}_{\mrm{red}}(\mcal{M}) . \]

\noindent 2.\
Multiplication induces a homomorphism
$\tilde{\mcal{A}}^{\bdot}_{X} \otimes_{\mcal{O}_{X}} \mcal{E}
\ar \tilde{\mcal{A}}^{\bdot}_{X}(\mcal{E})$.
Choose a {\em local algebraic frame}
$\bsym{f} = (f_{1}, \ldots, f_{r})^{\mrm{t}}$ for $\mcal{E}$ on
a small open set $U$; i.e.\ an isomorphism
$\bsym{f} : \mcal{O}_{U}^{r} \iso \mcal{E}|_{U}$.
Then we see that
$\tilde{\mcal{A}}^{\bdot}_{X} \otimes_{\mcal{O}_{X}} \mcal{E}
\iso \tilde{\mcal{A}}^{\bdot}_{X}(\mcal{E})$.

\noindent 3.\
The DO $\mrm{d} : \ul{\mbb{A}}(\mcal{M}) \ar
\ul{\mbb{A}}(\mcal{N})$
respects the cosimplicial structure.
\end{proof}

In order to clarify the algebraic structure of
$\tilde{\mcal{A}}^{\bdot}_{X}$
we introduce the following local objects. Given a chain
$\xi = (x_{0}, \ldots, x_{l})$ in $X$ let
\begin{equation} \label{eqn2.2}
\tilde{\mcal{A}}^{p,q}_{\xi} := \Omega^{q}(\Delta^{l}_{\mbb{Q}})
\otimes_{\mbb{Q}} \Omega^{p}_{X / k, \xi} .
\end{equation}
As usual we set
$\mrm{D}' := (-1)^{q} \otimes \mrm{d}$,
$\mrm{D}'' := \mrm{d} \otimes 1$ and
$\mrm{D} := \mrm{D}' + \mrm{D}''$.
The DGA $(\tilde{\mcal{A}}^{\bdot}_{\xi}, \mrm{D})$
is generated (as a DGA) by
\begin{equation} \label{eqn2.5}
\tilde{\mcal{A}}^{0}_{\xi} =
\mcal{O}_{X, \xi}[t_{0}, \ldots, t_{l}] / (\sum t_{i} - 1) .
\end{equation}
When $X$ is smooth of dimension $n$ over $k$ near $x_{0}$ then
$\tilde{\mcal{A}}^{p,q}_{\xi}$ is free of
rank $\binom{n}{p} \binom{l}{q}$ over $\tilde{\mcal{A}}^{0}_{\xi}$.
Given a quasi-coherent $\mcal{O}_{X}$-module $\mcal{M}$ let
$\tilde{\mcal{A}}^{p,q}_{\xi}(\mcal{M}) :=
\tilde{\mcal{A}}^{p,q}_{\xi} \otimes_{\mcal{O}_{X}} \mcal{M}$.

\begin{lem} \label{lem2.2}
\begin{enumerate}
\item For any quasi-coherent $\mcal{O}_{X}$-module $\mcal{M}$ and
open set $U \subset X$ there are natural commutative diagrams
\[ \begin{CD}
\Gamma(U, \tilde{\mcal{A}}^{p,q}_{X}(\mcal{M})) @>{\int_{\Delta}}>>
\Gamma(U, \mcal{A}^{p,q}_{X}(\mcal{M})) \\
@V{\Phi^{p, q}_{\mcal{M}}}VV  @VVV \\
\prod_{\xi \in S(U)} \tilde{\mcal{A}}^{p,q}_{\xi}(\mcal{M})
@>{\int_{\Delta}}>>
\prod_{\xi \in S(U)^{\mrm{red}}_{q}}
(\Omega^{p}_{X / k, \xi} \otimes_{\mcal{O}_{X}} \mcal{M})[-q]
\end{CD}   \]
\item $\Phi_{\mcal{M}} := \sum \Phi^{p, q}_{\mcal{M}}$
is injective and commutes with $\mrm{D}'$ and $\mrm{D}''$.
\item $\Phi_{\mcal{O}_{X}}$ is a DGA homomorphism and
$\Phi_{\mcal{M}}$ is $\Gamma(U, \tilde{\mcal{A}}^{\bdot}_{X})$-linear.
\end{enumerate}
\end{lem}

\begin{proof}
This is immediate from Definition \ref{dfn1.1} and formula
(\ref{eqn1.6}).
\end{proof}

\begin{lem} \label{lem2.3}
Let $\mcal{M}$ be a quasi-coherent sheaf. The natural homomorphism
$\mcal{M} \ar \tilde{\mcal{A}}^{0}_{X}(\mcal{M})$
extends to an $\mcal{O}_{X}$-linear homomorphism
$\ul{\mbb{A}}^{0}(\mcal{M}) \ar
\tilde{\mcal{A}}^{0}_{X}(\mcal{M})$.
\end{lem}

\begin{proof}
Consider the $i$-th covertex map
$\sigma_{i} : [0] \ar [l]$, which is the simplex
$\sigma_{i} = (i) \in \Delta^{l}_{0} \cong
\opn{Hom}_{\Delta}([0], [l])$
(cf.\ Example \ref{exa1.5}). There is a corresponding homomorphism
$\sigma_{i} : \ul{\mbb{A}}^{0}(\mcal{M}) \ar
\ul{\mbb{A}}^{l}(\mcal{M})$.
Given a local section $u \in \ul{\mbb{A}}^{0}(\mcal{M})$,
send it to
$(u_{0}, u_{1}, \ldots) \in
\tilde{\mrm{N}}^{0} \ul{\mbb{A}}(\mcal{M}) =
\tilde{\mcal{A}}^{0}_{X}(\mcal{M})$,
where
$u_{l} := \sum_{i=0}^{l} t_{i} \otimes \sigma_{i}(u)$.
\end{proof}

Because of the functoriality of our constructions we have:

\begin{prop} \label{prop2.3}
Let $f : X \ar Y$ be a morphism of $k$-schemes. Then the pullback
homomorphism
$f^{*} : \Omega^{\bdot}_{Y / k} \ar f_{*} \Omega^{\bdot}_{X / k}$
extends to DGA homomorphisms
$f^{*} : \tilde{\mcal{A}}^{\bdot}_{Y} \ar
f_{*} \tilde{\mcal{A}}^{\bdot}_{X}$
and
$f^{*} : \mcal{A}^{\bdot}_{Y} \ar f_{*} \mcal{A}^{\bdot}_{X}$
giving a commutative diagram
\[ \begin{CD}
\Omega^{\bdot}_{Y / k} @>>> \tilde{\mcal{A}}^{\bdot}_{Y}
@>{\int_{\Delta}}>> \mcal{A}^{\bdot}_{Y} \\
@V{f^{*}}VV @V{f^{*}}VV @V{f^{*}}VV \\
f_{*} \Omega^{\bdot}_{X / k} @>>> f_{*} \tilde{\mcal{A}}^{\bdot}_{X}
@>{f_{*}(\int_{\Delta})}>> f_{*} \mcal{A}^{\bdot}_{X} .
\end{CD} \]
\end{prop}

\begin{rem}
One can show that
$(\tilde{\mcal{A}}^{0}_{X})^{\times} = \mcal{O}_{X}^{\times}$
(invertible elements). We leave this as an exercise to the interested
reader.
\end{rem}

\section{Adelic Chern-Weil Theory}

Let us quickly review the notion of a connection on a module. For a
full account see \cite{GH} \S 0.5 and 3.3, \cite{KL} Appendix B,
\cite{Ka}, \cite{Go} or \cite{Kr}. In this section $k$ is a
field of characteristic $0$. Suppose
$A = A^{0} \oplus A^{1} \oplus \cdots$
is an associative, commutative DG $k$-algebra (i.e.\
$a b = (-1)^{ij} b a$ for $a \in A^{i}$, $b \in A^{j}$), with
operator $\mrm{d}$. Given an $A^{0}$-module $M$, a connection on
$M$ is a $k$-linear map
$\nabla : M \ar A^{1} \otimes_{A^{0}} M$ satisfying the Leibniz rule
$\nabla(a m) = \mrm{d} a \otimes m + a \nabla m$, $a \in A^{0}$.
$\nabla$ extends uniquely to an operator
$\nabla : A \otimes_{A^{0}} M \ar A \otimes_{A^{0}} M$
of degree $1$ satisfying the graded Leibnitz rule.

The curvature of $\nabla$ is the operator
$R := \nabla^{2} : M \ar A^{2} \otimes_{A^{0}} M$, which is
$A^{0}$-linear. The connection is flat, or integrable, if $R = 0$.
If $B$ is another DGA and $A \ar B$ is a DGA
homomorphism, then by extension of scalars there is an induced
connection
$\nabla_{B} : B^{0} \otimes_{A^{0}} M \ar B^{1} \otimes_{A^{0}} M$
over $B^{0}$.

If $M$ is free of rank $r$, choose a frame
$\bsym{e} = (e_{1}, \ldots, e_{r})^{\mrm{t}} : (A^{0})^{r} \iso M$.
Notice that we write $\bsym{e}$ as a column. This gives a
connection matrix $\bsym{\theta} = (\theta_{i, j})$,
$\theta_{i, j} \in A^{1}$, determined by
$\nabla \bsym{e} = \bsym{\theta} \otimes \bsym{e}$
(i.e.\
$\nabla e_{i} = \sum_{j} \theta_{i,j} \otimes e_{j}$).
In this case
$R \in A^{2} \otimes_{A^{0}} \opn{End}(M)$, so we get a curvature
matrix $\bsym{\Theta} = (\Theta_{i,j})$ satisfying
$R = \sum_{i,j} \Theta_{i,j} \otimes (e_{i} \otimes e^{\vee}_{j})$.
Here $\bsym{e}^{\vee} := (e_{1}^{\vee}, \ldots, e_{r}^{\vee})$ is
the dual basis and $\Theta_{i,j} \in A^{2}$.
One has
$\bsym{\Theta} = \mrm{d} \bsym{\theta} -
\bsym{\theta} \wedge \bsym{\theta}$.
If $\bsym{f}$ is another basis of $M$, with transition matrix
$\bsym{g} = (g_{i,j})$, $\bsym{e} = \bsym{g} \cdot \bsym{f}$,
then the matrix of
$\nabla$ w.r.t.\ $\bsym{f}$ is
$\bsym{g}^{-1} \bsym{\theta} \bsym{g} -
\bsym{g}^{-1} \mrm{d} \bsym{g}$,
and the curvature matrix is
$\bsym{g}^{-1} \bsym{\Theta} \bsym{g}$.

\begin{exa}
The Levi-Civita connection on $M$ determined by $\bsym{e}$, namely
$\nabla = (\mrm{d}, \ldots, \mrm{d})$, has matrix $\bsym{\theta} = 0$
and so is integrable. In terms of another basis
$\bsym{f} = \bsym{g} \cdot \bsym{e}$
the matrix will be $- \bsym{g}^{-1} \mrm{d} \bsym{g}$.
\end{exa}

Denote by $\mrm{M}_{r}(k)$ the algebra of matrices over
the field $k$ and
$\mrm{M}_{r}(k)^{*} :=$ \linebreak
$\opn{Hom}_{k}(\mrm{M}_{r}(k), k)$.
Then the symmetric algebra
$\mrm{S}(\mrm{M}_{r}(k)^{*})$
is the algebra of polynomial functions on $\mrm{M}_{r}(k)$.
The algebra
$I_{r}(k) := \mrm{S}(\mrm{M}_{r}(k)^{*})^{\mrm{Gl}_{r}(k)}$
of conjugation-invariant functions is generated by the
elementary invariant polynomials
$P_{1} = \opn{tr}, \ldots, P_{r} = \opn{det}$, with $P_{i}$
homogeneous of degree $i$.

\begin{lem} \label{lem3.1}
Assume that $A^{1} = A^{0} \cdot \mrm{d} A^{0}$.
Given any matrix
$\bsym{\theta} \in \mrm{M}_{r}(A^{1})$
let
$\bsym{\Theta} := \mrm{d} \bsym{\theta} -
\bsym{\theta} \cdot \bsym{\theta}$.
Then for any $P \in I_{r}(k)$
one has $\mrm{d} P(\bsym{\Theta}) = 0$,
\end{lem}

\begin{proof}
By assumption we can write
$\theta_{i, j} = \sum_{l} b_{i, j, l} \mrm{d} a_{l}$
for suitable $a_{l}, b_{i, j, l} \in A^{0}$.
Let $A_{\mrm{u}}$ be the universal algebra for this problem:
$A_{\mrm{u}}^{0}$ is the polynomial algebra
$k[ \bsym{a}, \bsym{b}]$,
where
$\bsym{a} = \{a_{l}\}$
and
$\bsym{b} = \{b_{i, j, l}\}$
are finite sets of indeterminates;
$A_{\mrm{u}} = \Omega^{\bdot}_{A_{\mrm{u}}^{0} / k}$;
and
$\bsym{\theta}_{\mrm{u}} \in \mrm{M}_{r}(A_{\mrm{u}}^{1})$
is the obvious connection matrix.
The DG $k$-algebra homomorphism
$A_{\mrm{u}} \ar A$ sends
$\bsym{\theta}_{\mrm{u}} \mapsto \bsym{\theta}$,
and hence it suffice to prove the case
$A = A_{\mrm{u}}$ .

Write $X := \opn{Spec} A^{0}$, which is nothing but
affine space $\mbf{A}^{N}_{k}$ for some $N$. We want to show
that the form
$\mrm{d} P(\bsym{\Theta})  = 0 \in
\Gamma(X, \Omega^{\bdot}_{X / k})$.
For a closed point $x \in X$ the residue field $k(x)$ is a
finite separable extension of $k$. This implies that the unique
$k$-algebra lifting
$k(x) \ar \widehat{\mcal{O}}_{X, x} = \mcal{O}_{X, (x)}$
has the property that
$\mrm{d} : \mcal{O}_{X, (x)} \ar \Omega^{1}_{X / k, (x)}$
is $k(x)$-linear. Since $X$ is smooth we have
$\mcal{O}_{X, (x)} \cong k(x)[[ f_{1}, \ldots, f_{N} ]]$.
We see that the differential equation on page 401 of \cite{GH} can
be solved formally in $\mrm{M}_{r}(\mcal{O}_{X, (x)})$.
Then the proof of the lemma on page 403 of \cite{GH} shows that
$\mrm{d} P(\bsym{\Theta})_{(x)} \in
\mfrak{m}_{x} \cdot \Omega^{\bdot}_{X / k, (x)}$.
Since this is true for all closed points $x \in X$ and
$\Omega^{\bdot}_{X / k, (x)}$ is a free $\mcal{O}_{X}$-module,
it follows that
$\mrm{d} P(\bsym{\Theta}) = 0$.
\end{proof}

Let us now pass to schemes. Assume $X$ is a finite type $k$-scheme
(not necessarily smooth), and let $\mcal{E}$ be a locally free
$\mcal{O}_{X}$-module of rank $r$.
We shall be interested in the sheaf of commutative
DGAs $\tilde{\mcal{A}}^{\bdot}_{X}$ and the locally free
$\tilde{\mcal{A}}^{0}_{X}$-module
$\tilde{\mcal{A}}^{0}_{X}(\mcal{E}) \cong
\tilde{\mcal{A}}^{0}_{X} \otimes_{\mcal{O}_{X}} \mcal{E}$.

\begin{dfn} \label{dfn3.2}
An {\em adelic connection} on $\mcal{E}$ is a connection
\[ \nabla : \tilde{\mcal{A}}^{0}_{X}(\mcal{E}) \ar
\tilde{\mcal{A}}^{1}_{X}(\mcal{E}) \]
over the algebra $\tilde{\mcal{A}}^{0}_{X}$.
\end{dfn}

\begin{dfn} \label{dfn3.3}
The {\em adelic curvature form} associated to an adelic connection
$\nabla$ on $\mcal{E}$ is
\[ R := \nabla^{2} \in
\opn{Hom}_{\tilde{\mcal{A}}_{X}^{0}}
\left( \tilde{\mcal{A}}_{X}^{0}(\mcal{E}),
\tilde{\mcal{A}}_{X}^{2}(\mcal{E}) \right)
\cong \Gamma \left( X, \tilde{\mcal{A}}_{X}^{2} \otimes_{\mcal{O}_{X}}
\mcal{E}nd_{\mcal{O}_{X}} (\mcal{E}) \right) . \]
\end{dfn}

Suppose $\nabla$ is an adelic connection on $\mcal{E}$
and $P \in I_{r}(k)$. Since
$P : \mcal{E}nd(\mcal{E}) \ar \mcal{O}_{X}$
is well defined, we get an induced sheaf homomorphism
$P : \tilde{\mcal{A}}^{\bdot}_{X} \otimes_{\mcal{O}_{X}}
\mcal{E}nd(\mcal{E}) \ar \tilde{\mcal{A}}^{\bdot}_{X}$
In particular we have
$P(R) \in \Gamma(X, \tilde{\mcal{A}}^{\bdot}_{X})$.

\begin{lem}
$P(R)$ is closed, i.e\ $\mrm{D} P(R) = 0$.
\end{lem}

\begin{proof}
This can be checked locally on $X$, so let $U$ be an open set on
which $\mcal{E}$ admits an algebraic frame
$\bsym{f}$. This frame induces isomorphisms of sheaves
$\bsym{f} : (\tilde{\mcal{A}}^{p, q}_{U})^{r} \iso
\tilde{\mcal{A}}^{p, q}_{X}(\mcal{E})|_{U}$
for all $p, q$. If
$\bsym{\theta} \in \mrm{M}_{r}(\Gamma(U,
\tilde{\mcal{A}}^{1}_{X}))$
is the matrix of the connection
$\nabla : \Gamma(U, \tilde{\mcal{A}}^{0}_{X}(\mcal{E})) \ar
\Gamma(U, \tilde{\mcal{A}}^{1}_{X}(\mcal{E}))$
then
$\bsym{\Theta} = \mrm{D} \bsym{\theta} - \bsym{\theta} \cdot
\bsym{\theta} \in
\mrm{M}_{r}(\Gamma(U, \tilde{\mcal{A}}^{2}_{X}))$
is the matrix of $R$, and we must show that
$\mrm{D} P(\bsym{\Theta}) = 0$.

According to Lemma \ref{lem2.2},
$\Phi : \Gamma(U, \tilde{\mcal{A}}^{\bdot}_{X}) \ar
\prod_{\xi \in S(U)} \tilde{\mcal{A}}^{\bdot}_{\xi}$
is an injective DGA homomorphism. Thus letting
$\bsym{\Theta}_{\xi}$ be the $\xi$-component of $\bsym{\Theta}$, it
suffices to show that $\mrm{D} P(\bsym{\Theta}_{\xi}) = 0$
for all $\xi$. Since
$\tilde{\mcal{A}}^{1}_{\xi} = \tilde{\mcal{A}}^{0}_{\xi} \cdot
\mrm{D} \tilde{\mcal{A}}^{0}_{\xi}$
we are done by Lemma \ref{lem3.1}.
\end{proof}

Recall that given a morphism of schemes $f : X \ar Y$ there is a
natural homomorphism of DGAs
$f^{*} : \tilde{\mcal{A}}^{\bdot}_{Y} \ar
f_{*} \tilde{\mcal{A}}^{\bdot}_{X}$.

\begin{prop}  \label{prop3.2}
Suppose $f : X \ar Y$ is a morphism of schemes, $\mcal{E}$ a
locally free $\mcal{O}_{Y}$-module and $\nabla$ an adelic connection
on $\mcal{E}$. Then there is an induced adelic connection
$f^{*}(\nabla)$ on $f^{*} \mcal{E}$, and
\[ f^{*}(P(R_{\nabla})) = P(R_{f^{*} (\nabla)}) \in
\Gamma(X, \tilde{\mcal{A}}^{\bdot}_{X}) . \]
\end{prop}

\begin{proof}
By adjunction there are homomorphisms
\[ f^{-1} \mcal{E} \xrightarrow{f^{-1}(\nabla)}
f^{-1} \tilde{\mcal{A}}^{1}_{Y}(\mcal{E}) \ar
\tilde{\mcal{A}}^{1}_{X}(f^{*} \mcal{E}) \]
of sheaves on $X$. Now
$\tilde{\mcal{A}}^{0}_{X}(f^{*} \mcal{E}) =
\tilde{\mcal{A}}^{0}_{X} \otimes_{f^{-1} \mcal{O}_{Y}}
f^{-1} \mcal{E}$,
so by Leibnitz rule we get $f^{*}(\nabla)$.
\end{proof}

\begin{thm} \label{thm3.2}
Let $X$ be a finite type
$k$-scheme and $\mcal{E}$ be a locally free $\mcal{O}_{X}$-module.
Choose an adelic connection $\nabla$ on $\mcal{E}$ and let $R$ be the
adelic curvature form. Then the $k$-algebra homomorphism
\tup{(}doubling degrees\tup{)}
\[ \begin{aligned}
w_{\mcal{E}} :
I_{r}(k) & \ar \mrm{H}^{\bdot} \Gamma(X, \tilde{\mcal{A}}^{\bdot}_{X})
\cong \mrm{H}^{\bdot}(X, \Omega^{\bdot}_{X / k}) \\
P & \mapsto [P(R)] .
\end{aligned} \]
is independent of the connection $\nabla$.
\end{thm}

We call $w_{\mcal{E}}$ the {\em adelic Chern-Weil homomorphism},
and we also write
$P(\mcal{E}) := w_{\mcal{E}}(P)$.

\begin{proof}
Suppose $\nabla'$ is another adelic connection, with curvature form
$R'$. We need to prove that
$[P(R)] = [P(R')] \in
\mrm{H}^{\bdot} \Gamma(X, \tilde{\mcal{A}}^{\bdot}_{X})$.

Consider the scheme
$Y := X \times \Delta^{1}_{\mbb{Q}}$, with projection morphisms
$s = s^{0} : Y \ar X$ and two sections
$\partial^{0}, \partial^{1} : X \ar Y$
(cf.\ Example \ref{exa1.3} for the notation).
Since
$s_{*} \Omega^{\bdot}_{Y / k} \cong
\Omega^{\bdot}_{X / k} \otimes_{\mbb{Q}}
\Omega^{\bdot}(\Delta^{1}_{\mbb{Q}})$
and
$\mbb{Q} \ar \Omega^{\bdot}(\Delta^{1}_{\mbb{Q}})$
is a quasi-isomorphism (Poincar\'{e} Lemma), we see that
$s^{*} : \Omega^{\bdot}_{X / k} \ar
s_{*} \Omega^{\bdot}_{Y / k}$
is a quasi-isomorphism of sheaves on $X$. Because $s$ is an
affine morphism the sheaves $\Omega^{p}_{Y / k}$ are acyclic for
$s_{*}$, and it follows that
$s_{*} \Omega^{\bdot}_{Y / k} \ar s_{*} \mcal{A}^{\bdot}_{Y}$
is a quasi-isomorphism. We conclude (cf.\ Proposition \ref{prop2.3})
that
$s^{*} : \mcal{A}^{\bdot}_{X} \ar s_{*} \mcal{A}^{\bdot}_{Y}$
is a quasi-isomorphism. Passing to global cohomology we also get
$\mrm{H}^{\bdot}(s^{*}) :
\mrm{H}^{\bdot} \Gamma(X, \mcal{A}^{\bdot}_{X}) \iso
\mrm{H}^{\bdot} \Gamma(Y, \mcal{A}^{\bdot}_{Y})$.
Therefore
\begin{equation} \label{eqn3.6}
\mrm{H}^{\bdot}(\partial^{0 *}) =
\mrm{H}^{\bdot}(\partial^{1 *}) :
\mrm{H}^{\bdot} \Gamma(Y, \tilde{\mcal{A}}^{\bdot}_{Y}) \iso
\mrm{H}^{\bdot} \Gamma(X, \tilde{\mcal{A}}^{\bdot}_{X})
\end{equation}
with inverse $\mrm{H}^{\bdot}(s^{*})$.

Let $\mcal{E}_{Y} := s^{*} \mcal{E}$, with two induced adelic
connections $s^{*} \nabla$ and $s^{*} \nabla'$.
Define the mixed adelic connection
\[ \nabla_{Y} := t_{0} s^{*} \nabla + t_{1} s^{*} \nabla' \]
on $\mcal{E}_{Y}$, with curvature $R_{Y}$. Now
$\partial^{0 *}(t_{0}) = 0$, so
\[ \partial^{0 *} \nabla_{Y} = \partial^{0 *}(t_{0} s^{*} \nabla)
+ \partial^{0 *}(t_{1} s^{*} \nabla') = \nabla'  \]
as connections on $\mcal{E}$.
Therefore
$\partial^{0 *}(P(R_{Y})) = P(R')$ and likewise
$\partial^{1 *}(P(R_{Y})) = P(R)$.
Finally use (\ref{eqn3.6}).
\end{proof}

Next we show how to construct adelic connections.

Recall that to every chain $\xi = (x_{0}, \ldots, x_{l})$
of length $l$ there is attached a DGA
\[ \tilde{\mcal{A}}^{\bdot}_{\xi} =
\Omega^{\bdot}(\Delta^{l}_{\mbb{Q}}) \otimes_{\mbb{Q}}
\Omega^{\bdot}_{X/k, \xi} \]
(cf.\ formula (\ref{eqn2.2})).
Set
$\tilde{\mcal{A}}^{i}_{\xi}(\mcal{E}) :=
\tilde{\mcal{A}}^{i}_{\xi} \otimes_{\mcal{O}_{X}} \mcal{E}$,
so $\tilde{\mcal{A}}^{0}_{\xi}(\mcal{E})$ is a free
$\tilde{\mcal{A}}^{0}_{\xi}$-module of rank $r$.
If $l = 0$ and $\xi = (x)$ then
$\tilde{\mcal{A}}^{\bdot}_{(x)} = \Omega^{\bdot}_{X/k, (x)}$
and
$\tilde{\mcal{A}}^{0}_{(x)} = \mcal{O}_{X, (x)} =
\widehat{\mcal{O}}_{X,x}$,
the complete local ring. For $0 \leq i \leq l$ there is a DGA
homomorphism
\[  \tilde{\mcal{A}}^{\bdot}_{(x_{i})} =
\Omega^{\bdot}_{X/k, (x_{i})} \xrightarrow{\sigma_{i}}
\Omega^{\bdot}_{X/k, \xi} \subset \tilde{\mcal{A}}^{\bdot}_{\xi} \]
(cf.\ proof of Lemma \ref{lem2.3}).

Suppose we are given a set
$\{ \nabla_{(x)} \}_{x \in X}$
where for each point $x$
\begin{equation} \label{eqn3.9}
\nabla_{(x)} : \mcal{E}_{(x)} \ar \Omega^{1}_{X / k, (x)}
\otimes_{\mcal{O}_{X, (x)}} \mcal{E}_{(x)}
\end{equation}
is a connection over $\mcal{O}_{X, (x)}$. Since
$\tilde{\mcal{A}}^{i}_{\xi}(\mcal{E}) \cong
\tilde{\mcal{A}}^{i}_{\xi} \otimes_{\mcal{O}_{X, (x)}}
\mcal{E}_{(x)}$,
each connection $\nabla_{(x_{i})}$ induces, by extension of scalars,
a connection
\[ \nabla_{\xi,i} :
\tilde{\mcal{A}}^{0}_{\xi}(\mcal{E}) \ar
\tilde{\mcal{A}}^{1}_{\xi}(\mcal{E}) \]
over the algebra $\tilde{\mcal{A}}^{0}_{\xi}$.
Define the ``mixed'' connection
\begin{equation} \label{eqn3.2}
\nabla_{\xi} := \sum_{i= 0}^{l} t_{i} \nabla_{\xi,i}:
\tilde{\mcal{A}}^{0}_{\xi}(\mcal{E}) \ar
\tilde{\mcal{A}}^{1}_{\xi}(\mcal{E}) .
\end{equation}

\begin{prop} \label{prop3.1}
Given a set of connections $\{ \nabla_{(x)} \}_{x \in X}$ as above,
there is a unique adelic connection
$\nabla$ on $\mcal{E}$, such that under the embedding
\[ \Phi_{\mcal{E}} :
\Gamma(U, \tilde{\mcal{A}}^{\bdot}_{X}(\mcal{E})) \subset
\prod_{\xi \in S(U)} \tilde{\mcal{A}}^{\bdot}_{\xi}(\mcal{E}) \]
of Lemma \tup{\ref{lem2.2}}, one has
$\nabla e = (\nabla_{\xi} e_{\xi})$
for every local section
$e = (e_{\xi}) \in$ \linebreak
$\Gamma(U, \tilde{\mcal{A}}^{0}_{X}(\mcal{E}))$.
Moreover,
$\nabla (\mcal{E}) \subset \tilde{\mcal{A}}^{1,0}_{X}(\mcal{E})$.
\end{prop}

\begin{proof}
The product
\[ \nabla := \prod_{\xi} \nabla_{\xi} :
\prod_{\xi \in S(X)} \tilde{\mcal{A}}^{0}_{\xi}(\mcal{E}) \ar
\prod_{\xi \in S(X)} \tilde{\mcal{A}}^{1}_{\xi}(\mcal{E}) \]
is a connection over the algebra
$\prod_{\xi} \tilde{\mcal{A}}^{0}_{\xi}$.
Since $\Phi_{\mcal{E}}$ is injective and $\Phi$ is a DGA homomorphism,
it suffices to show that
$\nabla e \in \tilde{\mcal{A}}^{1}_{X}(\mcal{E})$
for every local section
$e \in \tilde{\mcal{A}}^{0}_{X}(\mcal{E})$.

First consider a local section $e \in \mcal{E}$. For every point $x$,
\[ \nabla_{(x)} e \in \left( \Omega^{1}_{X/k}
\otimes_{\mcal{O}_{X}} \mcal{E} \right)_{(x)}  . \]
Therefore, writing
$\nabla_{l} := \prod_{\xi \in S(X)_{l}} \nabla_{\xi}$,
we see that
\[ \nabla_{0} e \in
\ul{\mbb{A}}^{0}(\mcal{O}_{X}) \otimes_{\mcal{O}_{X}}
\Omega^{1}_{X/k} \otimes_{\mcal{O}_{X}} \mcal{E} \cong
\ul{\mbb{A}}^{0}(\Omega^{1}_{X/k}
\otimes_{\mcal{O}_{X}} \mcal{E}) . \]
According to Lemma \ref{lem2.3} we get a section
\[ \bsym{\alpha} = (\alpha_{0} , \alpha_{1} , \ldots) \in
\tilde{\mcal{A}}^{0}_{X}(\Omega^{1}_{X/k}
\otimes_{\mcal{O}_{X}} \mcal{E}) \cong
\tilde{\mcal{A}}^{1,0}_{X}(\mcal{E})  \]
with
\[ \alpha_{l} =
\sum_{i=0}^{l} t_{i} \otimes \sigma_{i}(\nabla_{0} e) =
\nabla_{l} e , \]
so $\bsym{\alpha} = \nabla e$.

Finally, any section of
$\tilde{\mcal{A}}^{0}_{X}(\mcal{E}) \cong
\tilde{\mcal{A}}^{0}_{X} \otimes_{\mcal{O}_{X}} \mcal{E}$
is locally a sum of tensors
$a \otimes e$ with $a \in \tilde{\mcal{A}}^{0}_{X}$ and
$e \in \mcal{E}$, so by the Leibniz rule
\[ \nabla(a \otimes e) = \mrm{D} a \otimes e + a \nabla e \in
\tilde{\mcal{A}}^{1}_{X}(\mcal{E}) . \]
\end{proof}

Observe that relative to a local algebraic frame
$\bsym{f}$ for $\mcal{E}$,
the matrix of a  connection $\nabla$ as in the proposition
has entries in $\tilde{\mcal{A}}^{1, 0}_{X}$.

A {\em global adelic frame} for $\mcal{E}$ is a family
$\bsym{e} = \{ \bsym{e}_{(x)} \}_{x \in X}$,
where for each $x \in X$,
$\bsym{e}_{(x)} : \mcal{O}_{X, (x)}^{r} \iso \mcal{E}_{(x)}$
is a frame. In other words this is an isomorphism
$\bsym{e} : \ul{\mbb{A}}^{0}(\mcal{O}_{X})^{r} \iso
\ul{\mbb{A}}^{0}(\mcal{E})$
of $\ul{\mbb{A}}^{0}(\mcal{O}_{X})$-modules.

The next corollary is inspired by the work of Parshin \cite{Pa}.

\begin{cor} \label{cor3.2}
A global adelic frame $\bsym{e}$ of $\mcal{E}$ determines an
adelic connection $\nabla$.
\end{cor}

\begin{proof}
The frame $\bsym{e}_{(x)}$ determines a Levi-Civita connection
$\nabla_{(x)}$ on $\mcal{E}_{(x)}$. Now use Proposition
\ref{prop3.1}.
\end{proof}

We call such an connection {\em pointwise trivial}. In Sections
5 and 7 we shall only work with pointwise trivial
connections.

Given a local section
$\alpha \in \tilde{\mcal{A}}_{X}^{\bdot}(\mcal{M})$
we write
$\alpha = \sum \alpha^{p, q}$
with
$\alpha^{p, q} \in \tilde{\mcal{A}}_{X}^{p, q}(\mcal{M})$.
For a chain $\xi$ we write $\alpha_{\xi}$ for the $\xi$ component of
$\Phi_{\mcal{M}}(\alpha)$ (see Lemma \ref{lem2.2}).

\begin{lem} \label{lem3.9}
Let $\nabla$ be the pointwise trivial connection on $\mcal{E}$
determined by an adelic frame $\bsym{e}$. Let
$\xi = (x_{0}, \ldots, x_{l})$ be a chain, and let $\bsym{f}$ be any
frame of $\mcal{E}_{\xi}$. Write
$\bsym{e}_{(x_{i})} = \bsym{g}_{i} \cdot \bsym{f}$
for matrices
$\bsym{g}_{i} \in \opn{Gl}_{r}(\mcal{O}_{X, \xi})$, $0 \leq i \leq l$.
Then:
\begin{enumerate}
\item The connection matrix of $\nabla_{\xi}$ w.r.t.\ the frame
$\bsym{f}$ is
$\bsym{\theta} = - \sum t_{i} \bsym{g}_{i}^{-1}
\mrm{d} \bsym{g}_{i}$.
\item Let $\bsym{\Theta}^{1,1}$ be the matrix of the curvature form
$R^{1, 1}_{\xi}$ w.r.t.\ the frame
$\bsym{f} \otimes \bsym{f}^{\vee}$.
Then
\[ \bsym{\Theta}^{1,1} = - \sum \mrm{d} t_{i}
\wedge \bsym{g}_{i}^{-1} \mrm{d} \bsym{g}_{i} . \]
\end{enumerate}
\end{lem}

\begin{proof}
Direct calculation.
\end{proof}

\begin{dfn} \label{dfn3.1}
The $i$-th Chern forms of $\mcal{E}$ with respect to the
adelic connection $\nabla$ are
\[ \begin{aligned}
\tilde{c}_{i}(\mcal{E}, \nabla) & :=
P_{i}(R) \in \Gamma(X, \tilde{\mcal{A}}_{X}^{2i}) \\
c_{i}(\mcal{E}, \nabla) & :=
\int_{\Delta} P_{i}(R) \in \Gamma(X, \mcal{A}_{X}^{2i}) .
\end{aligned} \]
\end{dfn}

Let $t$ be an indeterminate, and define
$P_{t} := \sum_{i = 1}^{r} P_{i} t^{i} \in I_{r}(k)[t]$.

\begin{prop}[Whitney Sum Formula] \label{prop3.5}
Let $X$ be a finite type $k$-scheme. \linebreak
Suppose
$0 \ar \mcal{E}' \ar \mcal{E} \ar \mcal{E}'' \ar 0$
is a short exact sequence of locally free $\mcal{O}_{X}$-modules. Then
there exist adelic connections $\nabla', \nabla, \nabla''$ on
$\mcal{E}', \mcal{E}, \mcal{E}''$ respectively, with
corresponding curvature forms $R', R, R''$, s.t.\
\[ P_{t}(R) = P_{t}(R') \cdot P_{t}(R'') \in
\Gamma(X, \tilde{\mcal{A}}^{\bdot}_{X})[t] . \]
\end{prop}

\begin{proof}
For any point $x \in X$ choose a splitting
$\sigma_{(x)} : \mcal{E}_{(x)}'' \ar \mcal{E}_{(x)}$
of the sequence, as modules over $\mcal{O}_{X, (x)}$. Also choose
frames $\bsym{e}_{(x)}'$, $\bsym{e}_{(x)}''$ for
$\mcal{E}_{(x)}'$, $\mcal{E}_{(x)}''$ respectively. Let
$\bsym{e}_{(x)} := (\bsym{e}_{(x)}', \sigma_{(x)}(\bsym{e}_{(x)}''))$
be the resulting frame of $\mcal{E}_{(x)}$. Use the global adelic
frame $\bsym{e} = \{ \bsym{e}_{(x)} \}$ to define an adelic
connection $\nabla$ on $\mcal{E}$; and likewise define $\nabla'$ and
$\nabla''$.

In order to check that $P_{t}(R) = P_{t}(R') \cdot P_{t}(R'')$
it suffices, according to Lemma \ref{lem2.2}, to look separately
at each chain $\xi = (x_{0}, \ldots, x_{q})$.
Let
$\bsym{g}_{i} \in \mrm{Gl}_{r}(\mcal{O}_{X, \xi})$
be the transition matrix
$\bsym{e}_{(x_{i})} = \bsym{g}_{i} \cdot \bsym{e}_{(x_{q})}$.
Because of our special choice of frames the initial segment of each
frame $\bsym{e}_{(x_{i})}$ is a frame for the submodule
$\mcal{E}_{\xi}' \subset \mcal{E}_{\xi}$,
which implies that
$\bsym{g}_{i} = \left( \begin{smallmatrix}
\bsym{g}_{i}' & * \\ 0 & \bsym{g}_{i}''
\end{smallmatrix} \right)$,
where $\bsym{g}_{i}', \bsym{g}_{i}''$ are the obvious transition
matrices.

Now with respect to the frame $\bsym{e}_{(x_{q})}$ of
$\mcal{E}_{\xi}$, the connection matrix of $\nabla_{(x_{i})}$ is
$\bsym{\theta}_{i} = - \bsym{g}_{i}^{-1} \mrm{d} \bsym{g}_{i}$,
so the matrices of
$\nabla_{\xi}$ and $R_{\xi}$ are
\begin{eqnarray*}
\bsym{\theta} & = & - (t_{0} \bsym{g}_{0}^{-1} \mrm{d} \bsym{g}_{0}
+ \cdots + t_{q-1} \bsym{g}_{q-1}^{-1} \mrm{d} \bsym{g}_{q-1}) \\
\bsym{\Theta} & = & \mrm{D} \bsym{\theta} -
\bsym{\theta} \wedge \bsym{\theta} .
\end{eqnarray*}
It follows that
$\bsym{\Theta} = \left( \begin{smallmatrix}
\bsym{\Theta}' & * \\ 0 & \bsym{\Theta}''
\end{smallmatrix} \right)$.
By linear algebra we conclude that
$P_{t}(\bsym{\Theta}) = P_{t}(\bsym{\Theta}') \cdot
P_{t}(\bsym{\Theta}'')$.
\end{proof}

\begin{prop} \label{prop3.6}
If $\mcal{E}$ is an invertible $\mcal{O}_{X}$-module
and $\nabla$ is an adelic connection on it, then the differential
logarithm
\[ \opn{dlog} : \opn{Pic} X = \mrm{H}^{1}(X, \mcal{O}^{*}_{X})
\ar \mrm{H}^{2}(X, \Omega^{\bdot}_{X / k}) \cong
\mrm{H}^{2}(X, \mcal{A}^{\bdot}_{X}) \]
sends
\[ \opn{dlog}([\mcal{E}]) = [c_{1}(\mcal{E}; \nabla)] . \]
\end{prop}

\begin{proof}
(Cf.\ \cite{HY} Proposition 2.6.) Suppose $\{ U_{i} \}$
is a finite open cover of $X$ s.t.\ $\mcal{E}|_{U_{i}}$ is trivial
with frame $e_{i}$. Let
$g_{i, j} \in \Gamma(U_{i} \cap U_{j}, \mcal{O}^{*}_{X})$
satisfy $e_{i} = g_{i, j} e_{j}$. Then the \v{C}ech cocycle
$\{ g_{i, j} \} \in C^{1}(\{ U_{i} \}; \mcal{O}^{*}_{X})$
represents $[\mcal{E}]$.

Choose a global adelic frame $\{ e_{(x)} \}$ for $\mcal{E}$ and let
$\nabla$ be the connection it determines. For a chain $(x, y)$ let
$g_{(x, y)} \in \mcal{O}^{*}_{X, (x, y)}$ satisfy
$e_{(x)} = g_{(x, y)} e_{(y)}$. By Lemma \ref{lem3.9} we see that
$c_{1}(\mcal{E}; \nabla) = \{ \opn{dlog} g_{(x, y)} \} \in
\Gamma(X, \mcal{A}^{1,1}_{X})$.

For any point $x \in U_{i}$ define
$g_{i, (x)} \in \mcal{O}^{*}_{X, (x)}$ in the obvious way. Then
$\{ \opn{dlog} g_{i, (x)} \}$ \linebreak
$\in C^{0}(\{ U_{i} \}; \mcal{A}^{1, 0}_{X})$,
and
\[ \mrm{D} \{ \opn{dlog} g_{i, (x)} \} =
\{ \opn{dlog} g_{i, j} \} - \{ \opn{dlog} g_{(x, y)} \} . \]
Since
$\Gamma(X, \mcal{A}^{\bdot}_{X}) \ar
C^{\bdot}(\{ U_{i} \}; \mcal{A}^{\bdot}_{X})$
is a quasi-isomorphism we are done.
\end{proof}

\begin{thm}  \label{thm3.4}
Suppose $X$ is smooth over $k$, so that
$\mrm{H}^{\bdot} \Gamma(X, \tilde{\mcal{A}}_{X}^{\bdot}) =
\mrm{H}^{\bdot}_{\mrm{DR}}(X)$. Then the Chern classes
\[ c_{i}(\mcal{E}) := [ \tilde{c}_{i}(\mcal{E}, \nabla) ] \in
\mrm{H}^{2i}_{\mrm{DR}}(X) \]
coincide with the usual ones.
\end{thm}

\begin{proof}
By Theorem \ref{thm3.2} and Propositions \ref{prop3.2},
\ref{prop3.5} and \ref{prop3.6} we see that the axioms of
Chern classes (cf.\ \cite{Ha2} Appendix A) are satisfied.
\end{proof}

\begin{exa}
Consider the projective line $\mbf{P} = \mbf{P}^{1}_{k}$
and the sheaf $\mcal{O}_{\mbf{P}}(1)$.
Let $v \in \Gamma(\mbf{P}, \mcal{O}_{\mbf{P}}(1))$
have a zero at the point $z$. Define a global adelic frame
$\{ e_{(x)} \}$ by
$e_{(x)} = v$ if $x \neq z$, and
$e_{(z)} = w$, any basis of $\mcal{O}_{\mbf{P}}(1)_{(z)}$.
So $v = a w$ for some regular parameter
$a \in \mcal{O}_{\mbf{P}, (z)}$.
The local components of the Chern form
$c_{1}(\mcal{O}_{\mbf{P}}(1); \nabla)$ are $0$
unless $\xi = (z_{0}, z)$ ($z_{0}$ is the generic point),
where we get
$c_{1}(\mcal{O}_{\mbf{P}}(1); \nabla)_{\xi} = a^{-1} \mrm{d} a$.
\end{exa}

An algebraic connection on $\mcal{E}$ is a connection
$\nabla : \mcal{E} \ar \mcal{E} \otimes_{\mcal{O}_{X}}
\Omega^{1}_{X / k}$. The connection $\nabla$ is trivial if
$(\mcal{E}, \nabla) \cong (\mcal{O}_{X}, \mrm{d})^{r}$.
$\nabla$ is generically trivial if it's trivial on a dense open set.
The next proposition explores the relation between adelic and
algebraic connections.

\begin{prop} \label{prop3.10}
Assume $X$ is smooth irreducible and $k$ is algebraically
closed.
Let $\nabla$ be an integrable adelic connection on $\mcal{E}$.
\begin{enumerate}
\item If
$\nabla(\mcal{E}) \subset \tilde{\mcal{A}}^{1, 0}_{X}(\mcal{E})$
then $\nabla$ is algebraic.
\item If $\nabla$ is algebraic and generically trivial
then it is trivial.
\end{enumerate}
\end{prop}

\begin{proof}
1.\ By Lemma 2.15 part 1, with
$\mcal{M} = \Omega^{1}_{X / k} \otimes_{\mcal{O}_{X}} \mcal{E}$,
it suffices to prove that
$\mrm{D}'' \nabla(\mcal{E}) = 0$, which is
a local statement. So choose a local algebraic frame $\bsym{f}$ of
$\mcal{E}$ on some open set. Then we have a connection matrix
$\bsym{\theta}$ which is homogeneous of bidegree $(1, 0)$,
and by assumption the curvature matrix
$\bsym{\Theta} = \mrm{D} \bsym{\theta}
- \bsym{\theta} \cdot \bsym{\theta}$
is zero. But since
$\bsym{\Theta}^{1, 1} = \mrm{D}'' \bsym{\theta}$ we are done.

\medskip \noindent
2.\
The algebraic connection $\nabla$ extends uniquely to an adelic
connection with the same name (by Proposition \ref{prop3.1}).

Let $x_{0}$ be the generic point, so by assumption we have a frame
$\bsym{e}_{(x_{0})}$ for $\mcal{E}_{(x_{0})}$ which trivializes
$\nabla$.
Now take any closed point $x_{1}$, so
$\mcal{O}_{X, (x_{1})} \cong k[[t_{1}, \ldots, t_{n}]]$.
It is well known that there is a frame
$\bsym{e}_{(x_{1})}$ for $\mcal{E}_{(x_{1})}$ which trivializes
$\nabla$ (cf.\ \cite{Ka}).
Consider the chain $\xi = (x_{0}, x_{1})$. W.r.t.\ the frame
$\bsym{e}_{(x_{0})}$, the connection matrix of $\nabla_{\xi}$ is
$\bsym{\theta}_{\xi} = -t_{1} \bsym{g}^{-1} \mrm{d} \bsym{g}$,
where
$\bsym{e}_{(x_{1})} = \bsym{g} \cdot \bsym{e}_{(x_{0})}$
and
$\bsym{g} \in \mrm{Gl}_{r}(\mcal{O}_{X, \xi})$.

Since $\nabla$ is integrable we get
\[ -\mrm{d} t_{1} \cdot \bsym{g}^{-1} \mrm{d} \bsym{g} =
\mrm{D}'' \bsym{\theta}_{\xi} = \bsym{\Theta}^{1, 1} = 0 . \]
We conclude that
\[ \mrm{d} \bsym{g} = 0 \in \mrm{M}_{r}(\Omega^{1}_{X / k, \xi}) . \]
But because $X$ is smooth and $k$ is algebraically closed, it follows
that
$\mrm{H}^{0} \Omega^{\bdot}_{X / k, \xi} \subset
\mrm{H}^{0} \Omega^{\bdot}_{X / k, \eta} = k$,
where $\eta$ is any maximal chain containing $\xi$. So in fact
$\bsym{g} \in \mrm{Gl}_{r}(k)$, and by faithful flatness we get
\[ \bsym{e}_{(x_{0})} =
\bsym{g}^{-1} \cdot \bsym{e}_{(x_{1})} \in
\mcal{E}_{(x_{1})} \cap \mcal{E}_{(x_{0})} = \mcal{E}_{x_{1}} . \]

By going over all closed points $x_{1}$ we see that
$\bsym{e}_{(x_{0})} \in \Gamma(X, \mcal{E})$,
which trivializes $\nabla$.
\end{proof}

There do however exist integrable adelic connections which are not
algebraic.

\begin{exa} \label{exa3.10}
Let $X$ be any scheme of positive dimension, and let
$\tilde{a} \in \Gamma(X, \tilde{\mcal{A}}^{0}_{X})$
be any element satisfying
$\mrm{D}'' \mrm{D} \tilde{a} \neq 0$.
For instance take a fixed point $x_{0}$ and an element
$a_{(x_{0})} \in \mcal{O}_{X, (x_{0})}$
s.t.\ $\mrm{d} a_{(x_{0})} \neq 0$. For any $x \neq x_{0}$ set
$a_{(x)} := 0 \in \mcal{O}_{X, (x)}$.
Then
$\{ a_{(x)} \} \in \Gamma(X, \ul{\mbb{A}}^{0}(\mcal{O}_{X}))$,
and by Lemma \ref{lem2.3} we get
$\tilde{a} \in \Gamma(X, \tilde{\mcal{A}}^{0}_{X})$.
Now
$\mrm{D}'' \mrm{D} \tilde{a} = \mrm{D}'' \mrm{D}' \tilde{a}$,
and clearly $\mrm{D}' \tilde{a}$ is not algebraic.

Take $\mcal{E} = \mcal{O}_{X}^{2}$.
The matrix
\[ \bsym{e} =
\begin{pmatrix}
1 & \tilde{a} \\ 0 & 1 \end{pmatrix}
\in \mrm{M}_{2}(\Gamma(X, \tilde{\mcal{A}}^{0}_{X})) \]
is invertible, and we consider $\bsym{e}$ as a frame for
$\tilde{\mcal{A}}^{0}_{X}(\mcal{E})$.
Define $\nabla$ to be the Levi-Civita connection for $\bsym{e}$.
If we now consider the algebraic frame
$\bsym{f} = \left( \begin{smallmatrix}
1 & 0 \\ 0 & 1 \end{smallmatrix} \right)$
of $\mcal{E}$, then the connection matrix w.r.t.\ $\bsym{f}$ is
\[ -\bsym{e}^{-1} \cdot \mrm{D} \bsym{e} =
- \begin{pmatrix}
0 & \mrm{D} \tilde{a} \\ 0 & 0 \end{pmatrix}  .\]
So $\mrm{D}'' \nabla(\mcal{E}) \neq 0$ and hence
$\nabla$ is not algebraic.
\end{exa}

\begin{question} \label{que3.1}
Does there exist an adelic connection
$\nabla$ with curvature form $R$ homogeneous of bidegree $(1,1)$?
\end{question}

\begin{rem}
Theorem \ref{thm3.2} works just as well for a relative situation:
$Y$ is a finite type $k$-scheme and $f : X \ar Y$ is a finite
type morphism. Then we can define
$\mcal{A}^{\bdot}_{X / Y} :=
\mrm{N} \ul{\mbb{A}}(\Omega^{\bdot}_{X / Y})$
and likewise $\tilde{\mcal{A}}^{\bdot}_{X / Y}$.
There are relative adelic connections on any locally free
$\mcal{O}_{X}$-module $\mcal{E}$, and there is a Chern-Weil
homomorphism
\[ w_{\mcal{E}} : I_{r}(k) \ar
\mrm{H}^{\bdot} f_{*} \mcal{A}^{\bdot}_{X / Y} \cong
\mrm{H}^{\bdot} \mrm{R} f_{*} \Omega^{\bdot}_{X / Y} . \]
\end{rem}

\begin{rem}
In \cite{Du} a very similar construction is carried out to calculate
characteristic classes of principal $G$-bundles, for a Lie group $G$.
These classes are in the cohomology of the classifying space $BG$,
which coincides with the cohomology of the simplicial manifold $NG$.
\end{rem}

\begin{rem}
Suppose $X$ is any finite type scheme over $k$. Then
$R \in \Gamma(X, \tilde{\mcal{A}}_{X}^{\bdot})$
is nilpotent and we may define
\[ \tilde{\mrm{ch}}(\mcal{E}; \nabla) := \opn{tr} \opn{exp} R \in
\Gamma(X, \tilde{\mcal{A}}_{X}^{\bdot}) . \]
Using the idea of the proof of Proposition \ref{prop3.5} one can
show that given a bounded complex
$\mcal{E}_{\bdot}$ of locally free sheaves, which is acyclic on an
open set $U$, it is possible to find
connections $\nabla_{i}$ on $\mcal{E}_{i}$ s.t.\
$\sum_{i} (-1)^{i} \tilde{\mrm{ch}}(\mcal{E}_{i}; \nabla_{i}) = 0$
on $U$. In particular when $U = X$ we get a ring homomorphism
\[ \mrm{ch} : \mrm{K}^{0}(X) \ar
\mrm{H}^{\bdot} \Gamma(X, \tilde{\mcal{A}}_{X}^{\bdot}) \cong
\mrm{H}^{\bdot}(X, \Omega^{\bdot}_{X/k}) , \]
the {\em Chern character}. When $X$ is smooth this
is the usual Chern character into
$\mrm{H}^{\bdot}_{\opn{DR}}(X)$.
\end{rem}

\section{Secondary Characteristic Classes}

Let $k$ be a field of characteristic $0$.
In \cite{BE}, Bloch and Esnault show that given a locally free
sheaf $\mcal{E}$ on a smooth $k$-scheme $X$, an algebraic connection
$\nabla : \mcal{E} \ar \Omega^{1}_{X / k}
\otimes_{\mcal{O}_{X}} \mcal{E}$
and an invariant polynomial $P \in I_{r}(k)$ of degree $m \geq 2$,
there is a Chern-Simons class
\[ \mrm{T} P(\mcal{E}; \nabla) \in
\Gamma \left(X, \Omega^{2m -1}_{X / k} /
\mrm{d}(\Omega^{2m - 2}_{X / k}) \right) \]
satisfying
\[ \mrm{d} \mrm{T} P(\mcal{E}; \nabla) =
P(\mcal{E}) \in \mrm{H}^{2m}_{\mrm{DR}}(X) . \]
$\mrm{T} P(\mcal{E}; \nabla)$ is called the secondary, or
Chern-Simons, characteristic class. The notation we use is taken from
\cite{Es}; the original notation in \cite{BE} is
$w_{m}(\mcal{E}, \nabla, P)$.

Such algebraic connections exist when $X$ is affine. However
the authors of \cite{BE} point out that any quasi-projective scheme
$X$ admits a vector bundle whose total space $X'$ is affine, and then
$\mrm{H}^{\bdot}_{\mrm{DR}}(X) \ar \mrm{H}^{\bdot}_{\mrm{DR}}(X')$
is an isomorphism.

In Section 3 we proved that adelic connections always exist.
In this section we define adelic Chern-Simons classes, which are
global sections of sheaves on $X$ itself:

\begin{thm} \label{thm4.1}
Suppose $X$ is a smooth $k$-scheme, $\mcal{E}$ a locally free
$\mcal{O}_{X}$-module, $\nabla$ an adelic connection on $\mcal{E}$
and $P \in I_{r}(k)$ homogeneous of degree $m \geq 2$.
Then there is a class
\[ \mrm{T} P(\mcal{E}; \nabla) \in
\Gamma \left( X, \mcal{A}^{2m-1}_{X} /
\mrm{D}(\mcal{A}^{2m-2}_{X}) \right) \]
satisfying
\[ \mrm{D} [\mrm{T} P(\mcal{E}; \nabla)] = P(\mcal{E}) \in
\mrm{H}^{2m}_{\mrm{DR}}(X)  \]
and commuting with pullbacks for morphisms of schemes $X' \ar X$.
\end{thm}

The proof is later in this section, after some preparation.

According to \cite{BE} Theorem 2.2.1, for any commutative
$k$-algebra $B$, invariant polynomial
$P \in I_{r}(k)$ homogeneous of degree $m$ and matrix
$\bsym{\theta} \in \mrm{M}_{r}(\Omega^{1}_{B / k})$,
there is a differential form
$\mrm{T} P(\bsym{\theta}) \in \Omega^{2m-1}_{B / k}$.
$\mrm{T} P(\bsym{\theta})$ is functorial on $k$-algebras, and
satisfies
\begin{equation} \label{eqn4.1}
\mrm{d} \mrm{T} P(\bsym{\theta}) = P(\bsym{\Theta}) \in
\Omega^{2m}_{B / k} ,
\end{equation}
where
$\bsym{\Theta} = \mrm{d} \bsym{\theta} -
\bsym{\theta} \cdot \bsym{\theta}$.

We shall need a slight generalization of \cite{BE}
Proposition 2.2.2. Consider $\opn{M}_{r}$ and $\opn{GL}_{r}$
as schemes over $k$. There is a universal invertible matrix
\[ \bsym{g} = \bsym{g}_{\mrm{u}} \in
\Gamma (\mrm{GL}_{r}, \mrm{GL}_{r}(\mcal{O}_{\mrm{GL}_{r}})) . \]
For an integer $N$ there is a universal connection matrix
\[ \bsym{\theta}  = \bsym{\theta}_{\mrm{u}}
\in \Gamma(Y, \mrm{M}_{r}(\Omega^{1}_{Y / k})) , \]
where
$Y := \mrm{M}_{r} \times \mbf{A}^{N} =
\opn{Spec} [\bsym{a}_{\mrm{u}}, \bsym{b}_{\mrm{u}}]$
for a collection of indeterminates
$\bsym{a}_{\mrm{u}} = \{ a_{p} \}$
and
$\bsym{b}_{\mrm{u}} = \{ b_{i, j, p} \}$,
$1\leq p \leq N$, $1 \leq i, j \leq r$.
The matrix is of course
$\bsym{\theta}_{\mrm{u}} = (\theta_{i, j})$
with
$\theta_{i, j} = \sum_{p} b_{i, j, p} \mrm{d} a_{p}$.
Then we get by pullback matrices $\bsym{g}$ and $\bsym{\theta}$
on $\mrm{GL}_{r} \times Y$.

\begin{lem}
Given an invariant polynomial $P$ there is an open cover
$\mrm{GL}_{r} = \bigcup U_{i}$
and forms
$\beta_{i} = \beta_{\mrm{u}, i} \in \Gamma(U_{i} \times Y,
\Omega^{2m-2}_{U_{i} \times Y})$
s.t.\
\[ \left( \mrm{T} P(\bsym{\theta}) -
\mrm{T} P(\mrm{d} \bsym{g} \cdot \bsym{g}^{-1} +
\bsym{g} \cdot \bsym{\theta} \cdot \bsym{g}^{-1} ) \right)|_{
U_{i} \times Y} = \mrm{d} \beta_{i} . \]
\end{lem}

\begin{proof}
Write
\[ \alpha = \alpha_{\mrm{u}} := \mrm{T}P(\bsym{\theta}) -
\mrm{T}P(\mrm{d} \bsym{g} \cdot \bsym{g}^{-1} +
\bsym{g} \cdot \bsym{\theta} \cdot \bsym{g}^{-1})
\in \Gamma(\mrm{GL}_{r} \times Y,
\Omega^{2m - 1}_{\mrm{GL}_{r} \times Y / k}) . \]
It's known that $\mrm{d} \alpha = 0$. Let
$s : \mrm{GL}_{r} \ar \mrm{GL}_{r} \times Y$
correspond to any $k$-rational point of $Y$. Choose a point
$x \in \mrm{GL}_{r}$. By \cite{BE} Proposition 2.2.2 there is an
affine open neighborhood $V$ of $s(x)$ in $\mrm{GL}_{r} \times Y$
and a form
$\beta' \in \Gamma(V, \Omega^{2m-2}_{\mrm{GL}_{r} \times Y / k})$,
s.t.\
$\alpha|_{V} = \mrm{d} \beta'$.
Define $U := s^{-1}(V)$, so
$s^{*}(\alpha)|_{U} = \mrm{d}\, s^{*}(\beta') \in
\Gamma(U, \Omega^{2m-2}_{\mrm{GL}_{r} / k})$.
Since
$\mrm{H}(s^{*}) :
\mrm{H}^{\bdot}_{\mrm{DR}}(U \times Y) \ar
\mrm{H}^{\bdot}_{\mrm{DR}}(U)$
is an isomorphism, it follows that there is some
$\beta \in \Gamma(U \times Y,
\Omega^{2m-2}_{U \times Y / k})$
with
$\alpha|_{U \times Y} = \mrm{d} \beta$.
\end{proof}

\begin{proof}[Proof of the Theorem]
Say $\opn{dim} X = n$.
Let $U$ be a sufficiently small affine open set of $X$ s.t.\
$\mrm{d} a_{1}, \ldots, \mrm{d} a_{n}$ is an algebraic frame of
$\Omega^{1}_{X / k}$, for some
$a_{1}, \ldots, a_{n} \in \Gamma(U, \mcal{O}_{X})$;
and there a local algebraic frame $\bsym{f}$ for $\mcal{E}$ on $U$.

We get an induced isomorphism
$\bsym{f} : (\tilde{\mcal{A}}_{U}^{0})^{r} \iso
\tilde{\mcal{A}}^{0}_{X}(\mcal{E})|_{U}$.
Let
$\bsym{\theta} = (\theta_{i,j}) \in
\mrm{M}_{r}(\Gamma(U, \tilde{\mcal{A}}^{1}_{U}))$
be the connection matrix of $\nabla$ with respect to $\bsym{f}$.
Define the commutative DGAs
\[ A_{l} := \Omega^{\bdot}(\Delta^{l}_{\mbb{Q}}) \otimes_{\mbb{Q}}
\Gamma(U, \ul{\mbb{A}}^{l}(\Omega^{\bdot}_{X / k})) . \]
Then by definition \ref{dfn1.1},
$\bsym{\theta} = (\bsym{\theta}_{1}, \bsym{\theta}_{2}, \ldots)$
where
$\bsym{\theta}_{l} \in \mrm{M}_{r}(A_{l}^{1})$,
and the various matrices $\bsym{\theta}_{l}$ have to satisfy
certain simplicial compatibility conditions.

Fix an index $l$. We have
$A_{l}^{0} = \mcal{O}(\Delta^{l}_{\mbb{Q}}) \otimes_{\mbb{Q}}
\Gamma(U, \ul{\mbb{A}}^{l}(\mcal{O}_{X}))$
which contains $\Gamma(U, \mcal{O}_{X})$.
Thus we may uniquely write
$(\theta_{i,j})_{l} = \sum_{p = 1}^{n} b_{i,j,p} \mrm{d} a_{p} +
\sum_{p = 1}^{l} b_{i,j,p + n} \mrm{d} t_{p}$,
with $b_{i,j,p} \in A_{l}^{0}$.
It follows that for $N = n + l$
and $Y = \mrm{M}_{r} \times \mbf{A}^{N}$
there is a unique $k$-algebra homomorphism
$\phi_{l} : \Gamma(Y, \mcal{O}_{Y}) \ar A^{0}_{l}$,
with
$\phi_{l}(\bsym{a}_{\mrm{u}}, \bsym{b}_{\mrm{u}}) =
(\bsym{a}, \bsym{t}, \bsym{b})$.
This extends to a DGA homomorphism
$\phi_{l} : \Gamma(Y, \Omega^{\bdot}_{Y / k}) \ar A_{l}$,
and sends the universal connection $\bsym{\theta}_{\mrm{u}}$ to
$\bsym{\theta}_{l}$. Define
$\mrm{T}P(\bsym{\theta}_{l}) :=
\phi_{l}(\mrm{T}P(\bsym{\theta}_{\mrm{u}})) \in A_{l}$.

Because the homomorphisms $\phi_{l}$ are completely determined by the
matrices $\bsym{\theta}_{l}$ it follows that the forms
$\mrm{T}P(\bsym{\theta}_{l})$ satisfy the simplicial compatibilities.
So there is an adelic form
$\mrm{T}P(\bsym{\theta}) \in
\Gamma(U, \tilde{\mcal{A}}^{2m -1 }_{X})$.

Now let $\bsym{f}' = \bsym{g} \cdot \bsym{f}$ be another local
algebraic frame for $\mcal{E}$ on $U$, with
$\bsym{g} \in \Gamma(U, \mrm{GL}_{r}(\mcal{O}_{X}))$.
Fix $l$ as before and write
\[ \alpha_{l} := \mrm{T} P(\bsym{\theta}_{l}) -
\mrm{T} P(\mrm{d} \bsym{g} \cdot \bsym{g}^{-1} +
\bsym{g} \cdot \bsym{\theta}_{l} \cdot \bsym{g}^{-1})
\in \Gamma(U, \tilde{\mcal{A}}^{2m - 1}_{X}) . \]
Let $h : U \ar \mrm{GL}_{r}$ be the scheme morphism s.t.\
$h^{*} (\bsym{g}_{\mrm{u}}) = \bsym{g}$.
The $k$-algebra homomorphism
\[ \psi_{l} = h^{*} \otimes \phi_{l} :
\Gamma(\mrm{GL}_{r} \times Y, \mcal{O}_{\mrm{GL}_{r} \times Y})
\ar A^{0}_{l} \]
extends to a DGA homomorphism and
$\alpha_{l} = \psi_{l}(\alpha_{\mrm{u}})$,
where $\alpha_{\mrm{u}}$ is the obvious universal form.
By the lemma, for every $i$ there is a form
\[ \beta_{i, l} := \psi_{l}(\beta_{\mrm{u}, i}) \in
\Omega^{\bdot}(\Delta^{l}_{\mbb{Q}}) \otimes_{\mbb{Q}}
\Gamma(h^{-1}(U), \ul{\mbb{A}}^{l}(\Omega^{\bdot}_{X / k})) \]
of degree $2m -2$.

Since we are not making choices to define
$\beta_{i, l}$ it follows that the simplicial compatibilities hold,
and so we obtain an adele
$\beta_{i} \in \Gamma(h^{-1}(U), \tilde{\mcal{A}}^{2m-2}_{X})$,
which evidently satisfies
\[ \alpha|_{h^{-1}(U_{i})} = \mrm{D} \beta_{i}
\in \Gamma(h^{-1}(U_{i}), \tilde{\mcal{A}}^{2m-1}_{X}) . \]

This means that the element
\[ \mrm{T} \tilde{P}(\mcal{E}; \nabla) := \mrm{T} P(\bsym{\theta})
\in \Gamma \left( U, \tilde{\mcal{A}}^{2m-1}_{X}  /
\mrm{D}(\tilde{\mcal{A}}^{2m-1}_{X})
\right) \]
is independent of the local algebraic frame $\bsym{f}$,
and therefore glues to a global section on $X$.
Finally set
$\mrm{T} P(\mcal{E}; \nabla) := \int_{\triangle}
\mrm{T} \tilde{P}(\mcal{E}; \nabla)$.
\end{proof}

Some of the deeper results of \cite{BE} deal with integrable
algebraic connections. Denote by
$\mcal{H}^{i}_{\mrm{DR}}$ the sheafified De Rham cohomology, i.e.\
the sheaf associated to the presheaf
$U \mapsto \mrm{H}^{i}_{\mrm{DR}}(U)$.
Then
\[ \begin{aligned}
\mcal{H}^{i}_{\mrm{DR}} & = \mrm{H}^{i} \Omega^{\bdot}_{X / k} =
\opn{Ker} \left(
\frac{\Omega^{i}_{X / k} }{ \mrm{d}(\Omega^{i - 1}_{X / k}) }
\xrightarrow{\mrm{d}} \Omega^{i + 1}_{X / k} \right) \\
& \cong \mrm{H}^{i} \mcal{A}^{\bdot}_{X} =
\opn{Ker} \left(
\frac{ \mcal{A}^{i}_{X} }{ \mrm{D}( \mcal{A}^{i- 1}_{X}) }
\xrightarrow{\mrm{D}} \mcal{A}^{i + 1}_{X} \right) .
\end{aligned} \]
Because of formula (\ref{eqn4.1}), we get an adelic generalization of
\cite{BE} Proposition 2.3.2:

\begin{prop} \label{prop4.1}
If the adelic connection $\nabla$ is integrable then
$\mrm{T} P(\mcal{E}; \nabla) \in$ \linebreak
$\Gamma(X, \mcal{H}^{2m - 1}_{\mrm{DR}})$.
\end{prop}

The next question is an extension of Basic Question 0.3.1
of \cite{BE}.

\begin{question}
Are the classes $\mrm{T} P (\mcal{E}; \nabla)$ all zero for an
integrable adelic connection $\nabla$?
\end{question}

\section{The Bott Residue Formula}

Let $X$ be smooth $n$-dimensional projective variety over the field
$k$ ($\opn{char} k = 0 $). Suppose $\mcal{E}$ is a locally free
$\mcal{O}_{X}$-module of rank $r$ and
$P \in I_{r}(k)$ is a homogeneous polynomial of degree $n$.
The problem is to calculate the Chern number
$\int_{X} P(\mcal{E}) \in k$,
where
$\int_{X} : \mrm{H}^{2n}_{\mrm{DR}}(X) \ar k$
is the nondegenerate map defined e.g.\ in \cite{Ha1}.

Assume $v \in \Gamma(X, \mcal{T}_{X})$ is a vector field which acts
on $\mcal{E}$. By this we mean there is a DO
$\Lambda : \mcal{E} \ar \Omega^{1}_{X / k}
\otimes_{\mcal{O}_{X}} \mcal{E}$
satisfying
$\Lambda(a e) = v(a) \otimes e + a \otimes \Lambda e$
for local section $a \in \mcal{O}_{X}$ and $e \in \mcal{E}$.
Furthermore assume the zero scheme $Z$ of $v$ is finite
(but not necessarily reduced). Then we
shall define a local invariant $P(v, \mcal{E}, z) \in k$ for
every zero $z \in Z$, explicitly in terms of local coordinates,
in equation (\ref{eqn5.1}). Our result is:

\begin{thm}[Bott Residue Formula] \label{thm5.1}
\[ \int_{X} P(\mcal{E}) = \sum_{z \in Z} P(v, \mcal{E}, z) . \]
\end{thm}

The proof appears later in this section, after some preparation.
It follows the original proof of Bott \cite{Bo2}, but
using algebraic residues and adeles instead of complex-analytic
methods. This is made possible by Proposition \ref{prop2.4} and
Theorem \ref{thm3.4}.
We show that a good choice of adelic connection $\nabla$ on
$\mcal{E}$ enables one to localize the integral to the zero locus
$Z$. This is quite distinct from the proof of the Bott Residue
Formula in \cite{CL}, where classes in Hodge cohomology
$\mrm{H}^{q}(X, \Omega^{p}_{X / k})$ are considered, and integration
is done using Grothendieck's global duality theory.

Let us first recall the local cohomology residue map
\[ \opn{Res}_{\mcal{O}_{X, (z)} / k} :
\mrm{H}^{n}_{z} (\Omega^{n}_{X / k}) \ar k \]
of \cite{Li} and \cite{HK}.
Choose local coordinates $f_{1}, \ldots, f_{n}$ at $z$, so
$\mcal{O}_{X, (z)} \cong$ \linebreak
$k(z)[[ f_{1}, \ldots, f_{n}]]$.
Local cohomology classes are represented by generalized fractions.
Given
$a = \sum_{\bsym{i}} b_{\bsym{i}} \bsym{f}^{\bsym{i}} \in
\mcal{O}_{X, (z)}$
where $\bsym{i} = (i_{1}, \ldots, i_{n})$,
$b_{\bsym{i}} \in k(z)$ and
$\bsym{f}^{\bsym{i}} = f_{1}^{i_{1}} \cdots f_{n}^{i_{n}}$,
the residue is
\begin{equation} \label{eqn5.7}
\opn{Res}_{\mcal{O}_{X, (z)} / k}
\gfrac{a \cdot  \mrm{d} f_{1} \wedge \cdots \wedge \mrm{d} f_{n}}{
f_{1}^{i_{1}}, \ldots, f_{n}^{i_{n}}} =
\opn{tr}_{k(z) / k}(b_{i_{1} - 1, \ldots, i_{n} - 1})
\in k .
\end{equation}

Let $a_{1}, \ldots, a_{n} \in \mcal{O}_{X}$ be the unique
local sections near $z$ satisfying
$v = \sum a_{i} \frac{\partial}{\partial f_{i}}$.
Then by definition
\[ \mcal{O}_{Z, z} \cong k(z)[[ f_{1}, \ldots, f_{n}]] /
(a_{1}, \ldots, a_{n}) . \]
The DO $\Lambda$ restricts to an $\mcal{O}_{Z}$-linear
endomorphism $\Lambda |_{Z}$ of $\mcal{O}_{Z}
\otimes_{\mcal{O}_{X}} \mcal{E}$,
giving an element
$P(\Lambda |_{Z}) \in \mcal{O}_{Z, z}$.
Choose any lifting $P'$ of $P(\Lambda |_{Z})$ to
$\mcal{O}_{X, (z)}$, and define
\begin{equation} \label{eqn5.1}
P(v, \mcal{E}, z) := (-1)^{\binom{n + 1}{2}}
\opn{Res}_{\mcal{O}_{X, (z)} / k}
\gfrac{P' \cdot \mrm{d} f_{1} \wedge \cdots \wedge
\mrm{d} f_{n}}{ a_{1}, \ldots, a_{n}} .
\end{equation}
The calculation of (\ref{eqn5.1}), given the $a_{i}$ and $P'$,
is quite easy: first express these elements as power series in
$\bsym{f}$. The rules for manipulating generalized fractions
are the same as for ordinary fractions, so the denominator can
be brought to be $\bsym{f}^{\bsym{i}}$. Now use (\ref{eqn5.7}).

\begin{exa} \label{exa5.1}
Let $X := \mbf{P}^{1}_{k}$ and $\mcal{E} := \mcal{O}_{X}(1)$. Let
$\opn{Spec} k[f] \subset X$ be the complement of
one point (infinity), and let $z$ be the origin (i.e.\ $f(z) = 0$).
We embed $\mcal{E}$ in the function field $k(X)$ as the
subsheaf of functions with at most one pole at $z$. Now
$v := f^{2} \frac{\partial}{\partial f} \in \Gamma(X, \mcal{T}_{X})$
is a global vector field, and we see that its action on $k(X)$
preserves $\mcal{E}$. So the theorem applies with $\Lambda = v$.
Here is the calculation:
the zero scheme of $v$ is $Z = \opn{Spec} k[f] / (a)$ with
$a = f^{2}$. Since
$\Lambda(f^{-1}) = f^{2} \frac{\partial f^{-2}}{\partial f} = -1$
we see that
$P(\Lambda) = -f$, and so
\[ P(v, \mcal{E}, z) := (-1)^{3}
\opn{Res}_{k[[f]]  / k} \gfrac{-f \mrm{d} f}{f^{2}} = 1 , \]
as should be.
\end{exa}

\begin{exa} \label{exa5.2}
If $z$ is a simple zero (that is to say $\mcal{O}_{Z, z} = k(z)$)
we recover the familiar formula of Bott.
Denote by $\mrm{ad}(v)$ the adjoint action of $v$
on $\mcal{T}_{X}$. Then $\mrm{ad}(v)|_{z}$ is an invertible
$k(z)$-linear endomorphism of $\mcal{T}_{X}|_{k(z)}$. Its matrix
w.r.t.\ the frame
$( \frac{\partial}{\partial f_{1}}, \ldots,
\frac{\partial}{\partial f_{n}})^{\mrm{t}}$
is
$- \left( \frac{\partial a_{i}}{\partial f_{j}} \right)^{\mrm{t}}$,
and (\ref{eqn5.1}) becomes
\[ P(v, \mcal{E}, z) = (-1)^{\binom{n}{2}}
\opn{tr}_{k(z) / k} \left( \frac{ P(\Lambda |_{k(z)})}{
\opn{det}(\opn{ad}(v) |_{k(z)})} \right) .  \]

In the previous example we could have chosen
$v := f \frac{\partial}{\partial f}$. This has $2$ simple zeroes:
$z$, where $a = f$ and $P(\Lambda) = -1$, so
$P(v, \mcal{E}, z) = 1$; and infinity, where $P(\Lambda) = 0$.
\end{exa}

Let us start our proofs by showing that the local invariant is indeed
independent of choices.

\begin{lem} \label{lem5.1}
$P(v, \mcal{E}, z)$ is independent of the coordinate system
$f_{1}, \ldots, f_{n}$ and the lifting $P'$.
\end{lem}

\begin{proof}
Let $g_{1}, \ldots, g_{n}$ be another system of coordinates,
and write
$v = \sum b_{i} \frac{\partial}{\partial g_{i}}$.
Then we get
$a_{i} = \sum \frac{\partial f_{i}}{\partial g_{j}} b_{j}$.
The formulas for changing numerator and denominator in generalized
fractions imply that the value of (\ref{eqn5.1}) remains the same when
computed relative to $g_{1}, \ldots, g_{n}$.
\end{proof}

\begin{rem}
According to \cite{HK} Theorems 2.3 and 2.4 one has
$P(v, \mcal{E}, z) =$ \linebreak
$(-1)^{\binom{n}{2}}
\tau^{f}_{a} P(\Lambda|_{Z})$,
where
$\tau^{f}_{a} : \mcal{O}_{Z, z} \ar k$
is the trace of Scheja-Storch \cite{SS}.
\end{rem}

\begin{lem} \label{lem5.4}
There exists an open subset $U \subset X$ containing $Z$, and sections
$f_{1}, \ldots, f_{n} \in \Gamma(U, \mcal{O}_{X})$,
such that the corresponding morphism
$U \ar \mbf{A}^{n}_{k}$
is unramified \textup{(}and even \'{e}tale\textup{)}, and the fiber
over the origin is the reduced scheme $Z_{\mrm{red}}$. Thus
$\mcal{T}_{X}|_{U}$ has a frame
$( \frac{\partial}{\partial f_{1}}, \ldots,
\frac{\partial}{\partial f_{n}})^{\mrm{t}}$.
Moreover, we can choose $U$ s.t.\ there is a frame
$(e_{1}, \ldots, e_{r})^{\mrm{t}}$ for $\mcal{E}|_{U}$.
\end{lem}

\begin{proof}
As $X$ is projective and $Z$ finite we can certainly find an affine
open set $U = \opn{Spec} R$ containing $Z$.
For each point $z \in Z$ we can find sections
$f_{1, z}, \ldots, f_{n, z} \in R$
and
$e_{1, z}, \ldots, e_{r, z} \in \Gamma(U, \mcal{E})$
which satisfy the requirements at $z$.
Choose a ``partition of unity of $U$ to
order $1$ near $Z$'', i.e.\ a set of functions
$\{ \epsilon_{z} \}_{z \in Z} \subset R$ representing the
idempotents of $R / \sum_{z} \mfrak{m}_{z}^{2}$. Then define
$f_{i} := \sum_{z} \epsilon_{z} f_{i, z}$
and
$e_{i} := \sum_{z} \epsilon_{z} e_{i, z}$,
and shrink $U$ sufficiently.
\end{proof}

From here we continue along the lines of \cite{Bo2}, but of course
we use adeles instead of smooth functions.
The sheaf $\tilde{\mcal{A}}_{X}^{p, q}$ plays the role of the
sheaf of smooth $(p, q)$ forms on a complex manifold.
The operator $\mrm{D}''$ behaves like the anti-holomorphic
derivative $\bar{\partial}$; specifically $\mrm{D}'' \alpha = 0$ for
any $\alpha \in \Omega^{\bdot}_{X / k}$.

Fix an open set $U$ and sections $f_{1}, \ldots, f_{n}$ like in Lemma
\ref{lem5.4}. Then we get an algebraic frame
$(\frac{\partial}{\partial f_{1}}, \ldots,
\frac{\partial}{\partial f_{n}})^{\mrm{t}}$
of $\mcal{T}_{X}|_{U}$, and we can write the
vector field
$v = \sum a_{i} \frac{\partial}{\partial f_{i}}$
with $a_{i} \in \Gamma(U, \mcal{O}_{X})$.
Choose a global adelic frame $\{ \bsym{e}_{(x)} \}_{x \in X}$
for $\mcal{E}$ as follows:
\begin{equation} \label{eqn5.2}
\begin{array}{ll}
\bsym{e}_{(x)} =
(e_{1}, \ldots, e_{r})^{\mrm{t}} & \text{ if } x \in U \\
\bsym{e}_{(x)} = \text{ arbitrary } & \text{ if }  x \notin U .
\end{array}
\end{equation}
Then we get a family of connections
$\{ \nabla_{(x)} \}_{x \in X}$, and a global connection
$\nabla : \tilde{\mcal{A}}_{X}^{0}(\mcal{E}) \ar
\tilde{\mcal{A}}_{X}^{1}(\mcal{E})$
over the algebra $\tilde{\mcal{A}}_{X}^{0}$.
The curvature form
$R \in \tilde{\mcal{A}}_{X}^{2}(\mcal{E}nd(\mcal{E}))$
decomposes into homogeneous parts
$R = R^{2,0} + R^{1,1}$. Since $\tilde{\mcal{A}}_{X}^{p, q} = 0$ for
$p > n$, we get $P(R) = P(R^{1, 1})$; we will work with $R^{1,1}$.

\begin{lem} \label{lem5.2}
Applying the $\mcal{O}_{X}$-linear homomorphism
$\mrm{D}'' : \tilde{\mcal{A}}_{X}^{p, q}(
\mcal{E}nd (\mcal{E})) \ar$ \linebreak
$\tilde{\mcal{A}}_{X}^{p, q+1}(\mcal{E}nd (\mcal{E}))$
one has
$\mrm{D}'' R^{1,1} = 0$.
\end{lem}

\begin{proof}
This is a local statement. Passing to matrices using a local algebraic
frame, it is enough to prove that
$\mrm{D}'' \bsym{\Theta}^{1,1} = 0$.
Now $\bsym{\Theta} = \mrm{D} \bsym{\theta} -
\bsym{\theta} \wedge \bsym{\theta}$,
and $\bsym{\theta} \in \mrm{M}_{r}(\tilde{\mcal{A}}_{X}^{1,0})$,
so
$\bsym{\Theta}^{1,1} = \mrm{D}'' \bsym{\theta}$. But
$(\mrm{D}'')^{2} = 0$.
\end{proof}

Denote the canonical pairing
$\mcal{T}_{X} \otimes_{\mcal{O}_{X}} \Omega^{1}_{X / k}
\ar \mcal{O}_{X}$
by $\langle -, - \rangle$.
It extends to a bilinear pairing
$\tilde{\mcal{A}}_{X}^{0}(\mcal{T}_{X}) \otimes_{\mcal{O}_{X}}
\tilde{\mcal{A}}_{X}^{0}(\Omega^{1}_{X/k}) \ar
\tilde{\mcal{A}}_{X}^{0}$.
For each point $x \in X$ we choose a form
$\pi_{(x)} \in \Omega^{1}_{X / k, (x)}$ as follows:
\begin{equation} \label{eqn5.4}
\parbox{10cm}{ \begin{enumerate}
\item If $x \in Z$ set $\pi_{(x)} := 0$.
\item If $x \in U - Z$, let $j$ be the first index s.t.\
$a_{j}(x) \neq 0$, and set
$\pi_{(x)} := a_{j}^{-1} \mrm{d} f_{j}$.
\item If $x \notin U$ take any form
$\pi_{(x)} \in \Omega^{1}_{X / k, (x)}$ satisfying
$\langle v, \pi_{(x)} \rangle = 1$.
\end{enumerate} }
\end{equation}
Together we get a global section
$\pi = \{ \pi_{(x)} \}_{x \in X} \in
\ul{\mbb{A}}^{0}(\Omega^{1}_{X / k})$, and
as indicated in Lemma \ref{lem2.3}, there is a corresponding global
section
$\pi \in \tilde{\mcal{A}}_{X}^{0}(\Omega^{1}_{X / k}) =
\tilde{\mcal{A}}_{X}^{1, 0}$.

\begin{lem}
Considering
$v \in \tilde{\mcal{A}}_{X}^{0}(\mcal{T}_{X})$, one has the identity
\begin{equation} \label{eqn5.5}
\langle v, \pi \rangle = 1 \in \tilde{\mcal{A}}_{X}^{0}
\text{ on } X - Z .
\end{equation}
\end{lem}

\begin{proof}
Use the embedding of Lemma \ref{lem2.2} to reduce formula
(\ref{eqn5.5}) to the local formula
$\langle v, \pi_{\xi} \rangle = \sum t_{i} = 1
\in \tilde{\mcal{A}}_{\xi}^{0}$.
\end{proof}

In Bott's language (see \cite{Bo2}) $\pi$ is a projector for $v$.

Let
$\iota_{v} = \langle v, - \rangle: \Omega^{1}_{X / k}
\ar \mcal{O}_{X}$
be the interior derivative, or contraction along $v$. It extends
to an $\mcal{O}_{X}$-linear operator of degree $-1$ on
$\Omega^{\bdot}_{X / k}$,
and hence to an $\tilde{\mcal{A}}_{X}^{0}$-linear
operator of bidegree $(-1,0)$ on $\tilde{\mcal{A}}_{X}^{\bdot}$,
which commutes (in the graded sense)
with $\mrm{D}''$ and satisfies $\iota_{v}^{2} = 0$.

\begin{lem} \label{lem5.6}
There exists a global section
$L \in \tilde{\mcal{A}}_{X}^{0}(\mcal{E}nd(\mcal{E}))$
satisfying
\[ \iota_{v} R^{1, 1} = \mrm{D}'' L \in
\tilde{\mcal{A}}_{X}^{0,1}(\mcal{E}nd(\mcal{E})) \]
and
\[ L|_{Z} = \Lambda|_{Z} \in \mcal{E}nd(\mcal{E}|_{Z}) . \]
\end{lem}

\begin{proof}
Using Lemma \ref{lem2.1}, define
\[ L := \Lambda - \iota_{v} \circ \nabla :
\tilde{\mcal{A}}_{X}^{0}(\mcal{E})
\ar \tilde{\mcal{A}}_{X}^{0}(\mcal{E}) . \]
This is an $\tilde{\mcal{A}}_{X}^{0}$-linear homomorphism.
Let us distinguish between
$\mrm{D}'' L$, which is the image of $L$ under
$\mrm{D}'' :
\tilde{\mcal{A}}_{X}^{0}(\mcal{E}nd_{\mcal{O}_{X}}(\mcal{E}))
\ar
\tilde{\mcal{A}}_{X}^{0, 1}(
\mcal{E}nd_{\mcal{O}_{X}}(\mcal{E}))$,
and
$\mrm{D}'' \circ L$, which is the composed operator
$\tilde{\mcal{A}}_{X}^{0}(\mcal{E}) \ar
\tilde{\mcal{A}}_{X}^{0, 1}(\mcal{E})$.
Both $\iota_{v} R^{1, 1}$ and $\mrm{D}'' L$ can be thought of as
$\mcal{O}_{X}$-linear homomorphisms
$\mcal{E} \ar \tilde{\mcal{A}}_{X}^{0,1}(\mcal{E})$.
Since $\mrm{D}'' (\mcal{E}) = 0$, one checks (using a local algebraic
frame) that
$\mrm{D}'' L = \mrm{D}'' \circ L$
on $\mcal{E}$. By the proof of Lemma \ref{lem5.2},
$\mrm{D}'' \circ \nabla  = R^{1, 1}$
as operators
$\mcal{E} \ar \tilde{\mcal{A}}_{X}^{1, 1}(\mcal{E})$.
Now
$\mrm{D}'' \circ \Lambda = \Lambda \circ \mrm{D}''$.
Therefore we get equalities
\[ \mrm{D}'' L = \mrm{D}'' \circ L - L \circ \mrm{D}'' =
- \mrm{D}'' \circ \iota_{v} \circ \nabla
= \iota_{v} \circ \mrm{D}'' \circ \nabla  =
\iota_{v} R^{1, 1}  \]
of maps $\mcal{E} \ar \tilde{\mcal{A}}_{X}^{0,1}(\mcal{E})$.
Finally the equality
$L|_{Z} = \Lambda|_{Z}$ follows from the vanishing
of $\iota_{v}$ on $Z$.
\end{proof}

Let $t$ be an indeterminate, and define
\[ \begin{aligned}
\eta & := P(L + t R^{1, 1}) \cdot \pi \cdot
(1 - t \mrm{D}'' \pi)^{-1}  \\[1mm]
& = P(L + t R^{1, 1}) \cdot \pi
\cdot (1 + t \mrm{D}'' \pi + (t \mrm{D}'' \pi)^{2} + \cdots)
\in \tilde{\mcal{A}}_{X}^{\bdot} \sqbr{t}
\end{aligned} \]
(note that $(\mrm{D}'' \pi)^{n+1} = 0$, so this makes sense).
Writing
$\eta = \sum_{i} \eta_{i} t^{i}$ we see that
$\eta_{i} \in \tilde{\mcal{A}}_{X}^{i+1, i}$.

\begin{lem} \label{lem5.5}
$\mrm{D}'' \eta_{n-1} + P(R^{1,1}) = 0$ on $X - Z$.
\end{lem}

\begin{proof}
Using the multilinear polarization $\tilde{P}$ of $P$, Lemma
\ref{lem5.6} and the fact that $\iota_{v} - t \mrm{D}''$ is
an odd derivation, one sees that
$(\iota_{v} - t \mrm{D}'') P(L + t R^{1, 1}) = 0$
(cf.\ \cite{Bo2}). Since
$\langle v, \pi \rangle = 1$ on $X - Z$ we get
$(\iota_{v} - t \mrm{D}'') \pi = (1 - t \mrm{D}'') \pi$,
$(\iota_{v} - t \mrm{D}'')(1 - t \mrm{D}'' \pi) = 0$,
and hence
$(\iota_{v} - t \mrm{D}'') \eta = P(L + t R^{1, 1})$
on $X - Z$. Finally consider the coefficient of $t^{n}$ in this
expression, noting that $\eta_{n} = 0$, being a section of
$\tilde{\mcal{A}}_{X}^{n+1, n}$.
\end{proof}

The proof of the next lemma is easy.

\begin{lem} \label{lem3.2}
Let $P(M_{1}, \ldots, M_{n})$ be a multilinear polynomial on
$\mrm{M}_{r}(k)$, invariant under permutations.
Let $A = \bigoplus A^{i}$ be a commutative DG $k$-algebra
and $A^{-} := \bigoplus A^{2i+1}$.
Let $\alpha_{1}, \ldots, \alpha_{n} \in A^{-}$ and
$M_{1}, \ldots, M_{n} \in \opn{M}_{r}(A^{-})$. Then
\begin{enumerate}
\item If $\alpha_{i} = \alpha_{j}$ or $M_{i} = M_{j}$ for two distinct
indices $i,j$ then
\[ P(\alpha_{1} M_{1}, \ldots, \alpha_{n} M_{n}) = 0 . \]
\item
\[ P \left( \sum_{i=1}^{n} \alpha_{i} M_{i}, \ldots,
\sum_{i=1}^{n} \alpha_{i} M_{i} \right) =
n! P(\alpha_{1} M_{1}, \ldots, \alpha_{n} M_{n})  . \]
\end{enumerate}
\end{lem}

\begin{proof}[Proof of Theorem \ref{thm5.1}]
By definition
$c_{i}(\mcal{E}) = [ \int_{\Delta} \tilde{c}_{i}(\mcal{E}; \nabla) ]
\in \mrm{H}^{2i}_{\mrm{DR}}(X)$.
It is known that $\tilde{\mcal{A}}^{\bdot}_{X}$ is a commutative DGA,
and that
$\mrm{H}(\int_{\Delta}) :
\mrm{H} \Gamma(X, \tilde{\mcal{A}}^{\bdot}_{X}) \ar
\mrm{H}^{\bdot}_{\mrm{DR}}(X)$
is an isomorphism of graded algebras (see Corollary \ref{cor2.2}).
Hence
\[ Q(c_{1}(\mcal{E}), \ldots, c_{r}(\mcal{E})) =
[ \int_{\Delta} Q(\tilde{c}_{1}(\mcal{E}; \nabla), \ldots,
\tilde{c}_{r}(\mcal{E}; \nabla)) ]  =
[ \int_{\Delta} P(R) ] . \]
As mentioned before,
$P(R) = P(R^{1,1}) \in \tilde{\mcal{A}}_{X}^{2n}$.
We must verify:
\[ \int_{X} \int_{\Delta} P(R^{1,1}) = \sum_{z \in Z}
P(v, \mcal{E}, z) . \]
Let
\[
\Xi := S(U)_{n}^{\mrm{red}} - S(U - Z)_{n}^{\mrm{red}} =
\{ (x_{0}, \ldots, x_{n}) \ |\ x_{n} \in Z \} . \]
We are given that $X$ is proper, so by \cite{Be} (or by the
Parshin-Lomadze Residue Theorem, \cite{Ye1} Theorem 4.2.15)
\[ \int_{X} \int_{\Delta} \mrm{D}'' \eta_{n-1} =
\int_{X} \mrm{D}'' \int_{\Delta} \eta_{n-1} = 0 . \]
Therefore by Lemma \ref{lem5.5}
\[ \begin{aligned}
\int_{X} \int_{\Delta} P(R^{1,1}) & =
\int_{X} \int_{\Delta} (P(R^{1,1}) + \mrm{D}'' \eta_{n-1}) \\[1mm]
& = \sum_{\xi \in \Xi} \opn{Res}_{\xi} \int_{\Delta}
(P(R^{1,1}) + \mrm{D}'' \eta_{n-1}) .
\end{aligned} \]

Let us look what happens on the open set $U$. By construction the
connection $\nabla$ is integrable there
(it is a Levi-Civita connection w.r.t.\ the algebraic frame $\ul{e}$),
so $R = 0$ ; hence $P(R^{1,1}) = 0$ and
$\mrm{D}'' \eta_{n-1} = P(L) (\mrm{D}'' \pi)^{n}$.
According to Lemma \ref{lem5.6},
$\mrm{D}'' L = 0$, so by Lemma \ref{lem2.1} one has
$L \in \mcal{E}nd(\mcal{E})$. Therefore $P(L) \in \mcal{O}_{X}$.
Since $\int_{\Delta}$ is $\mcal{O}_{X}$-linear we get
$\int_{\Delta} \mrm{D}'' \eta_{n-1} =
P(L) \int_{\Delta} (\mrm{D}'' \pi)^{n}$.
All the above is on $U$. The conclusion is:
\begin{equation} \label{eqn5.6}
\int_{X} \int_{\Delta} P(R^{1,1}) =
\sum_{\xi \in \Xi} \opn{Res}_{\xi} \left( P(L) \int_{\Delta}
(\mrm{D}'' \pi)^{n} \right) .
\end{equation}

Using the embedding of Lemma \ref{lem2.2}, for each $\xi \in \Xi$
one has $\pi_{\xi} = \sum t_{i} \pi_{(x_{i})}$, and therefore
$\mrm{D}'' \pi_{\xi} = \sum \mrm{d} t_{i} \wedge \pi_{(x_{i})}$.
Let
\[ \Xi_{a} := \{ \xi = (x_{0}, \ldots, x_{n})\ |\
a_{1}(x_{i}) = \cdots = a_{i}(x_{i}) = 0 \text{ for }
i = 1, \ldots, n \} . \]
If $\xi \notin \Xi_{a}$ then for at least one index $0 \leq i < n$,
$\pi_{(x_{i})} = \pi_{(x_{i+1})}$. So by Lemma \ref{lem3.2} we get
$(\mrm{D}'' \pi_{\xi})^{n} = 0$.

It remains to consider only $\xi \in \Xi_{a}$. Since
$\pi_{(x_{i})} = a_{i+1}^{-1} \mrm{d} f_{i+1}$
for $0 \leq i < n$, and $\pi_{(x_{n})} = 0$, it follows from
Lemma \ref{lem3.2} that
\[ (\mrm{D}'' \pi_{\xi})^{n} = n! (-1)^{\binom{n}{2}}
\mrm{d} t_{0} \wedge \cdots \wedge \mrm{d} t_{n-1}
\wedge
\frac{\mrm{d} f_{1} \wedge \cdots \wedge
\mrm{d} f_{n}}{a_{1} \cdots a_{n}} , \]
so
\[ \int_{\Delta} (\mrm{D}'' \pi_{\xi})^{n} =
(-1)^{\binom{n+1}{2}}
\frac{\mrm{d} f_{1} \wedge \cdots \wedge
\mrm{d} f_{n}}{a_{1} \cdots a_{n}} \in \Omega^{n}_{X/k, \xi} . \]
Finally, according to \cite{Hu2} Corollary 2.5 or \cite{SY}
Theorem 0.2.9, and our Lemma \ref{lem5.1},
and using the fact that $L|_{Z} = \Lambda|_{Z}$, it holds:
\[ \begin{split}
 \sum_{\xi \in \Xi_{a}} (-1)^{l} & \opn{Res}_{\xi}
\frac{ P(L) \mrm{d} f_{1} \wedge \cdots \wedge \mrm{d} f_{n}}{
a_{1} \cdots a_{n}} \\[1mm]
& = \sum_{z \in Z} (-1)^{l} \opn{Res}_{\mcal{O}_{X (z)}/k}
\gfrac{ P( L ) \mrm{d} f_{1} \wedge \cdots \wedge \mrm{d} f_{n} }{
a_{1}, \ldots, a_{n}} \\[1mm]
& = \sum_{z \in Z} P(v, \mcal{E}, z) .
\end{split} \]
\end{proof}

\begin{rem}
There is a sign error in \cite{Ye1} Section 2.4. Let $K$ and
$L =$ \linebreak
$K((t_{1}, \ldots, t_{n}))$ be topological local fields.
Since the residue map
$\opn{Res}_{L / K} : \Omega^{\bdot, \mrm{sep}}_{L/k} \ar
\Omega^{\bdot, \mrm{sep}}_{K/k}$
is an $\Omega^{\bdot, \mrm{sep}}_{K/k}$-linear map of degree $-n$,
it follows by transitivity that \linebreak
$\opn{Res}_{L / K}(t_{1}^{-1} \mrm{d} t_{1} \wedge \cdots \wedge
t_{n}^{-1} \mrm{d} t_{n})= 1$, not
$(-1)^{\binom{n}{2}}$. This error carried into \cite{SY}
Theorem 0.2.9.
\end{rem}

\begin{rem}
Consider an action of an algebraic torus $T$ on $X$ and $\mcal{E}$,
with positive dimensional fixed point locus.
A residue formula is known in this case (see \cite{Bo2} and
\cite{AB}), but we were unable to prove it using our adelic method.
The sticky part was finding an adelic
model for the equivariant cohomology $\mrm{H}^{\bdot}_{T}(X)$.
An attempt to use an adelic version of the Cartan-De Rham complex did
not succeed. Another unsuccessful try was to consider the
classifying space $BT$ as a cosimplicial scheme, and compute its
cohomology via adeles.

Recently Edidin and Graham defined equivariant Chow groups, and
pro\-ved the Bott formula in that context (see \cite{EG}).
The basic idea there
is to approximate the classifying space $BT$ by finite type schemes,
``up to a given codimension''.
This approach is suited for global constructions, but
again it did not help as far as adeles were concerned.
\end{rem}

\section{The Gauss-Bonnet Formula}

Let $k$ be a perfect field (of any characteristic) and $X$ a finite
type scheme, with structural morphism $\pi : X \ar \opn{Spec} k$.
According to Grothendieck Duality Theory \cite{RD} there is a functor
$\pi^{!} : \msf{D}^{+}_{\mrm{c}}(\msf{Mod}(k)) \ar
\msf{D}^{+}_{\mrm{c}}(\msf{Mod}(X))$
between derived categories, called the twisted inverse image. The
object $\pi^{!} k$ is a dualizing complex on $X$, and it has a
canonical representative, namely the residue complex
$\mcal{K}^{\bdot}_{X} := \mrm{E} \pi^{!} k$. Here $\mrm{E}$ is the
Cousin functor. As graded sheaf
$\mcal{K}^{-q}_{X} = \bigoplus_{\opn{dim} \overline{ \{ x \} } = q}
\mcal{K}_{X}(x)$,
where $\mcal{K}_{X}(x)$ is a quasi-coherent sheaf, constant with
support $\overline{ \{x\} }$, and as $\mcal{O}_{X, x}$-module it
is an injective hull of the residue field $k(x)$.

$\mcal{K}^{\bdot}_{X}$ enjoys some remarkable properties,
which are deduced from corresponding properties of $\pi^{!}$.
If $X$ is pure of dimension $n$ then there is a canonical
homomorphism
$C_{X} : \Omega^{n}_{X / k} \ar \mcal{K}^{-n}_{X}$,
which gives a quasi-isomorphism
$\Omega^{n}_{X / k}[n] \ar \mcal{K}^{\bdot}_{X}$
on the smooth locus of $X$. If $f : X \ar Y$ is a morphism of schemes
then there is a map of graded sheaves
$\opn{Tr}_{f} : f_{*} \mcal{K}^{\bdot}_{X} \ar \mcal{K}^{\bdot}_{Y}$,
which becomes a map of complexes when $f$ is proper.

For integers $p, q$ define
\[ \mcal{F}_{X}^{p, q} := \mcal{H}om_{X}
(\Omega^{-p}_{X / k}, \mcal{K}^{q}_{X}) . \]
Clearly $\mcal{F}^{\bdot}_{X}$ is a graded (left and right)
$\Omega^{\bdot}_{X / k}$-module. Moreover according to \cite{Ye3}
Theorem 4.1, or \cite{EZ}, $\mcal{F}^{\bdot}_{X}$ is in fact a DG
$\Omega^{\bdot}_{X / k}$-module. The difficult thing is to define the
dual operator
$\opn{Dual}(\mrm{d}) : \mcal{F}_{X}^{p, q} \ar
\mcal{F}_{X}^{p + 1, q}$,
which is a differential operator of order $1$.
$\mcal{F}^{\bdot}_{X}$ is called the {\em De Rham-residue complex}.

There is a special cocycle
$C_{X} \in \Gamma(X, \mcal{F}_{X}^{\bdot})$
called the {\em fundamental class}. When $X$ is integral of
dimension $n$ then $C_{X} \in \mcal{F}_{X}^{-n, -n}$
is the natural map
$\Omega^{n}_{X/k} \ar \mcal{K}^{-n}_{X} = k(X) \otimes_{\mcal{O}_{X}}
\Omega^{n}_{X/k}$.
For any closed subscheme $f: Z \ar X$, the trace map
$\mrm{Tr}_{f} : f_{*} \mcal{F}_{Z}^{\bdot} \ar
\mcal{F}_{X}^{\bdot}$
is injective, which allows us to write
$C_{Z} \in \mcal{F}_{X}^{\bdot}$.
If $Z_{1}, \ldots, Z_{r}$ are the irreducible components of $Z$
(with reduced subscheme structures) and
$z_{1}, \ldots, z_{r}$ are the generic points,
one has
$C_{Z} = \sum_{i} (\opn{length} \mcal{O}_{Z, z_{i}}) C_{Z_{i}}$.

When $X$ is pure of dimension $n$,
$C_{X}$ induces a map of complexes
$C_{X} : \Omega^{\bdot}_{X / k}[2n]$ \linebreak
$\ar \mcal{F}^{\bdot}_{X}$,
$\alpha \mapsto C_{X} \cdot \alpha$.
This is a quasi-isomorphism on the smooth locus of $X$ (\cite{Ye3}
Proposition 5.8). Hence when $X$ is smooth
$\mrm{H}^{-i} \Gamma(X, \mcal{F}^{\bdot}_{X}) \cong
\mrm{H}_{i}^{\mrm{DR}}(X)$,
De Rham homology.

Given a morphism $f : X \ar Y$, the trace
$\opn{Tr}_{f} : f_{*} \mcal{K}^{\bdot}_{X}
\ar \mcal{K}^{\bdot}_{Y}$
and
$f^{*} : \Omega^{\bdot}_{Y / k} \ar f_{*} \Omega^{\bdot}_{X / k}$
induce a map of graded sheaves
$\opn{Tr}_{f} : f_{*} \mcal{F}^{\bdot}_{X}
\ar \mcal{F}^{\bdot}_{Y}$.
Given an \'{e}tale morphism $g : U \ar X$ there is a homomorphism of
complexes
$\mrm{q}_{g} : \mcal{F}^{\bdot}_{X} \ar g_{*} \mcal{F}^{\bdot}_{U}$,
and $\mrm{q}_{g}(C_{X}) = C_{U}$.

The next theorem summarizes a few theorems in \cite{Ye5} about
the action of $\mcal{A}^{\bdot}_{X}$ on $\mcal{F}^{\bdot}_{X}$.

\begin{thm} \label{thm6.1}
Let $X$ be a finite type scheme over a perfect field $k$.
\begin{enumerate}
\item $\mcal{F}_{X}^{\bdot}$ is a right DG
$\mcal{A}^{\bdot}_{X}$-module, and the
multiplication extends the $\Omega^{\bdot}_{X / k}$-module structure.
\item If $f : X \ar Y$ is proper then
$\opn{Tr}_{f} : f_{*} \mcal{F}^{\bdot}_{X} \ar \mcal{F}^{\bdot}_{Y}$
is $\mcal{A}^{\bdot}_{Y}$-linear.
\item If $g : U \ar X$ is \'{e}tale then
$\mrm{q}_{g} : \mcal{F}^{\bdot}_{X} \ar g_{*} \mcal{F}^{\bdot}_{U}$
is  $\mcal{A}^{\bdot}_{X}$-linear.
\end{enumerate}
\end{thm}

Note that from part 1 it follows that if $X$ is smooth of dimension
$n$ then
$\mcal{A}^{\bdot}_{X}[2n] \ar \mcal{F}^{\bdot}_{X}$,
$\alpha \mapsto C_{X} \cdot \alpha$ is a quasi-isomorphism.

Let us say a few words about the multiplication
$\mcal{F}^{\bdot}_{X} \otimes \mcal{A}^{\bdot}_{X}
\ar \mcal{A}^{\bdot}_{X}$.
Since
$\mcal{A}^{\bdot}_{X} \cong
\ul{\mbb{A}}_{\mrm{red}}^{\bdot}(\mcal{O}_{X}) \otimes_{\mcal{O}_{X}}
\Omega^{\bdot}_{X / k}$
and
$\mcal{F}^{\bdot}_{X} \cong \mcal{H}om_{\mcal{O}_{X}}
(\Omega^{\bdot}_{X / k}, \mcal{K}^{\bdot}_{X})$,
it suffices to describe the product
$\mcal{K}^{\bdot}_{X} \otimes
\ul{\mbb{A}}_{\mrm{red}}^{\bdot}(\mcal{O}_{X}) \ar
\mcal{K}^{\bdot}_{X}$.
This requires the explicit construction
of $\mcal{K}_{X}^{\bdot}$ which we gave in \cite{Ye3}, and which
we quickly review below.

The construction
starts with the theory of {\em Beilinson completion algebras}
(BCAs) developed in \cite{Ye2}. A BCA $A$ is a semilocal $k$-algebra
with a topology and with valuations on its residue fields. Each
local factor of $A$ is a quotient of the ring of formal power series
$L((s_{n})) \cdots ((s_{1}))[[t_{1}, \ldots, t_{m}]]$,
where $L$ is a finitely generated extension field of $k$,
and $L((s_{n})) \cdots ((s_{1}))$
is the field of iterated Laurent series.
One considers two kinds of homomorphisms between BCAs:
morphisms $f : A \ar B$ and intensifications $u : A \ar \widehat{A}$.

Each BCA $A$ has a {\em dual module} $\mcal{K}(A)$, which is
functorial w.r.t.\ these homomorphisms; namely there are maps
$\opn{Tr}_{f} : \mcal{K}(B) \ar \mcal{K}(A)$
and
$\mrm{q}_{u} : \mcal{K}(A) \ar \mcal{K}(\widehat{A})$.
If $A$ is local with maximal ideal $\mfrak{m}$ and residue field
$K = A / \mfrak{m}$, then a choice of coefficient field
$\sigma : K \ar A$ determines an isomorphism
\[ \mcal{K}(A) \cong \opn{Hom}_{K}^{\mrm{cont}}(A,
\Omega^{n, \mrm{sep}}_{K / k}) , \]
where $\Omega^{\bdot , \mrm{sep}}_{K / k}$ is the separated algebra
of differentials on $K$ and
$n = \opn{rank}_{K} \Omega^{1, \mrm{sep}}_{K / k}$.
In particular, algebraically $\mcal{K}(A)$ is an injective hull of
$K$.

Suppose $\xi = (x, \ldots, y)$ is a saturated chain of points in
$X$ (i.e.\ immediate specializations). Then the Beilinson completion
$\mcal{O}_{X, \xi}$ is a BCA. The natural algebra homomorphisms
$\partial^{-} : \mcal{O}_{X, (x)} \ar \mcal{O}_{X, \xi}$ and
$\partial^{+} : \mcal{O}_{X, (y)} \ar \mcal{O}_{X, \xi}$
are an intensification and a morphism, respectively. So there are
homomorphisms on dual modules
$\mrm{q}_{\partial^{-}} : \mcal{K}(\mcal{O}_{X, (x)})
\ar \mcal{K}(\mcal{O}_{X, \xi})$
and
$\opn{Tr}_{\partial^{+}} : \mcal{K}(\mcal{O}_{X, \xi}) \ar
\mcal{K}(\mcal{O}_{X, (y)})$.
The composition
$\opn{Tr}_{\partial^{+}} \circ \opn{q}_{\partial^{-}}$ is denoted by
$\delta_{\xi}$. We regard
$\mcal{K}_{X}(x) := \mcal{K}(\mcal{O}_{X, (x)})$
as a quasi-coherent
$\mcal{O}_{X}$-module, constant on the closed set
$\overline{\{ x \}}$. Define
\begin{equation} \label{eqn1.1}
 \mcal{K}_{X}^{q} := \bigoplus_{\opn{dim} \overline{\{x\}} = -q}
\mcal{K}_{X}(x)
\end{equation}
and
\begin{equation} \label{eqn6.2}
\delta = (-1)^{q+1} \sum_{(x,y)} \delta_{(x,y)} :
\mcal{K}_{X}^{q} \ar \mcal{K}_{X}^{q + 1} .
\end{equation}
Then the pair $(\mcal{K}^{\bdot}_{X}, \delta)$ is the residue
complex of $X$. That is to say, there is a canonical isomorphism
$\mcal{K}^{\bdot}_{X} \cong \pi^{!} k$ in the derived
category $\msf{D}(\msf{Mod}(X))$ (see \cite{Ye3} Corollary 2.5).

Let $x$ be a point of dimension $q$ in $X$, and consider a local
section
$\phi_{x} \in \mcal{K}_{X}(x) \subset \mcal{K}_{X}^{-q}$.
Let $\xi = (x_{0}, \ldots, x_{q'})$ be any chain of length $q'$
in $X$, and let
$a_{\xi} \in \mcal{O}_{X, \xi}$.
Define
$\phi_{x} \cdot a_{\xi} \in \mcal{K}_{X}^{-q + q'}$
as follows.
If $x = x_{0}$ and $\xi$ is saturated, then
\[ \phi_{x} \cdot a_{\xi} :=
\opn{Tr}_{\partial^{+}}(a_{\xi} \cdot \mrm{q}_{\partial^{-}}
(\phi_{x})) \in \mcal{K}_{X}(x_{q'}) , \]
where the product
$a_{\xi} \cdot \mrm{q}_{\partial^{-}}(\phi_{x})$
is in $\mcal{K}(\mcal{O}_{X, \xi})$.
Otherwise set
$\phi_{x} \cdot a_{\xi} := 0$.

It turns out that for local sections
$\phi = (\phi_{x}) \in \mcal{K}^{-q}_{X}$
and
$a = (a_{\xi}) \in \ul{\mbb{A}}_{\mrm{red}}^{q'}(\mcal{O}_{X})$,
one has
$\phi_{x} \cdot a_{\xi} = 0$
for all but finitely many pairs $x, \xi$. Hence
\[ \phi \cdot a := \sum_{x, \xi} \phi_{x} \cdot a_{\xi} \in
\mcal{K}^{-q + q'}_{X} \]
is well defined, and this is the product we use.

\begin{exa} \label{exa6.1}
Suppose $X$ is integral of dimension $n$, $x_{0}$ is its generic
point and
$\phi \in \mcal{K}^{-n}_{X} = \mcal{K}_{X}(x_{0}) =
\Omega^{n}_{k(x_{0}) / k}$.
Consider a saturated chain
$\xi = (x_{0}, \ldots, x_{q})$
and an element
$a \in \mcal{O}_{X, \xi} = k(x_{0})_{\xi}$.
We want to see what is
$\psi := \phi \cdot a \in \mcal{K}_{X}(x_{q}) =
\mcal{K}(\mcal{O}_{X, (x_{q})})$.
Choose a coefficient field
$\sigma : k(x_{q}) \ar \mcal{O}_{X, (x_{q})}$,
so that
$\mcal{K}(\mcal{O}_{X, (x_{q})}) \cong
\opn{Hom}^{\mrm{cont}}_{k(x_{q})}(\mcal{O}_{X, (x_{q})},
\Omega^{n - q}_{k(x_{q}) / k})$.
It's known that $k(x_{0})_{\xi} = \prod L_{i}$, a finite product of
topological local fields (TLFs), and $\sigma : k(x_{q}) \ar L_{i}$
is a morphism of TLFs. Then for $b \in \mcal{O}_{X, (x_{q})}$
one has
\[ \psi(b) = \sum_{i} \opn{Res}_{L_{i} / k(x_{q})} (b a \phi) \in
\Omega^{n - q}_{k(x_{q}) / k} , \]
where $\opn{Res}_{L_{i} / k(x_{q})}$ is the the residue of \cite{Ye1}
Theorem 2.4.3, and the product
$b a \phi \in \Omega^{n, \mrm{sep}}_{k(x_{0}) / k}$.
\end{exa}

We can now state the main result of this section. From here to the
end of Section 6 we assume $\opn{char} k = 0$.

\begin{thm}[Gauss-Bonnet] \label{thm7.1}
Assume $\opn{char} k = 0$, and let $X$ be an integral,
$n$-dimensional, quasi-projective $k$-variety
\textup{(}not necessarily smooth\textup{)}.
Let $\mathcal{E}$ be a locally free $\mathcal{O}_{X}$-module of
rank $r$. Suppose $v \in \Gamma(X, \mathcal{E})$
is a regular section, with zero scheme $Z$.
Then there is an adelic connection $\nabla$ on $\mcal{E}$ satisfying
\[ C_{X} \cdot c_{r}(\mathcal{E}, \nabla) =
(-1)^{m}   C_{Z} \in \mathcal{F}_{X}^{-2(n - r)}  \]
with $m = nr + \binom{r+1}{2}$.
\end{thm}

Let $U$ be an open subset such that $\mcal{E}|_{U}$ is trivial
and $U$ meets each irreducible component of $Z$. Fix an
algebraic frame $(v_{1}, \ldots, v_{r})^{\mrm{t}}$ of
$\mcal{E}|_{U}$ and write
\begin{equation}
v = \sum_{i = 1}^{r} a_{i} v_{i},\ a_{i} \in \Gamma(U, \mcal{O}_{X}) .
\end{equation}

For each $x \in X$ choose a local frame $\bsym{e}_{(x)}$ of
$\mcal{E}_{(x)}$ as follows:
\begin{equation} \label{eqn7.2}
\parbox{10cm}{ \begin{enumerate}
\item If $x \notin Z \cup U$, take
$\bsym{e}_{(x)}  = (v, *, \ldots, *)^{\mrm{t}}$.
\item If $x \in U - Z$, there is some $0 \leq i < r$ such that
$a_{1}(x), \ldots, a_{i}(x) = 0$ but $a_{i+1}(x) \neq 0$. Then take\\
$\bsym{e}_{(x)}  = (v, v_{1}, \ldots, v_{i}, v_{i+2}, \ldots,
v_{r})^{\mrm{t}}$.
\item If $x \in Z \cap U$, take
$\bsym{e}_{(x)}  = (v_{1}, \ldots, v_{r})^{\mrm{t}}$.
\item If $x \in Z - U$, take $\bsym{e}_{(x)} $ arbitrary.
\end{enumerate} }
\end{equation}
Let $\nabla_{(x)}$ be the resulting Levi-Civita connection on
$\mcal{E}_{(x)}$, let
$\nabla : \tilde{\mcal{A}}_{X}^{0}(\mcal{E}) \ar
\tilde{\mcal{A}}_{X}^{1}(\mcal{E})$ be the induced adelic connection,
and let
$R \in \tilde{\mcal{A}}_{X}^{2}(\mcal{E}nd(\mcal{E}))$
be the curvature. We get a top Chern form
$P_{r}(R) = \opn{det} R \in \tilde{\mcal{A}}_{X}^{2r}$.
Under the embedding of DGAs
$\Gamma(X, \tilde{\mcal{A}}_{X}^{p,q}) \subset \prod_{\xi \in S(X)}
\tilde{\mcal{A}}_{\xi}^{p,q}$
of Lemma \ref{lem2.2} we write $R = (R_{\xi})$.

\begin{lem} \label{lem7.1}
Suppose $\xi = (x_{0}, \ldots, x_{q})$ is a saturated chain of
length $q$,
with $x_{0}$ the generic point of $X$, and either:
\textup{(i)}\ $q < r$;
\textup{(ii)}\ $x_{q} \notin Z$; or
\textup{(iii)}\ $q = r$ and $\bsym{e}_{x_{i}} = \bsym{e}_{x_{i+1}}$
for some $i$. Then $\int_{\Delta^{q}} \opn{det} R_{\xi} = 0$.
\end{lem}

\begin{proof}
Let
$\bsym{g}_{i} \in \opn{GL}_{r}(\mcal{O}_{X, \xi})$ be the transition
matrix
$\bsym{e}_{x_{i}} = \bsym{g}_{i} \cdot \bsym{e}_{x_{q}}$,
let $\bsym{\theta}_{i}$ be the connection matrix of $\nabla_{x_{i}}$
w.r.t.\ the frame $\bsym{e}_{x_{q}}$. Then the matrices
$\bsym{\theta}$, $\bsym{\Theta}$ of $\nabla_{\xi}$, $R_{\xi}$
are
\[ \begin{aligned}
\bsym{\theta} & =  - (t_{0} \bsym{g}_{0}^{-1} \mrm{d} \bsym{g}_{0}
+ \cdots + t_{q-1} \bsym{g}_{q-1}^{-1} \mrm{d} \bsym{g}_{q-1}) \\
\bsym{\Theta} & =  \mrm{D} \bsym{\theta} -
\bsym{\theta} \wedge \bsym{\theta} .
\end{aligned} \]
In cases (i) and (ii), all $x_{i} \notin Z$, so
\[ g_{i} =
\left( \begin{smallmatrix}
1 & 0 & \dots & 0 \\
\vdots & \vdots & \ddots & \vdots \\
* & * & \dots & * \\
* & * & \dots & *
\end{smallmatrix} \right)
\hspace{5mm}
\bsym{\theta}_{i}  =
\left( \begin{smallmatrix}
0 & 0 & \dots & 0 \\
\vdots & \vdots & \ddots & \vdots \\
* & * & \dots & * \\
* & * & \dots & *
\end{smallmatrix} \right)
. \]
Therefore $\bsym{\Theta}$ has a zero first row too and
$\opn{det} \bsym{\Theta} = 0$.

Now suppose $q = r$. Since
$\bsym{\Theta} \in \opn{M}_{r}(\tilde{\mcal{A}}_{\xi}^{1, 1}) \oplus
\opn{M}_{r}(\tilde{\mcal{A}}_{\xi}^{2, 0})$, from degree
considerations we conclude that
$\int_{\Delta^{r}} \opn{det} \bsym{\Theta} =
\int_{\Delta^{r}} \opn{det} (\bsym{\Theta}^{1,1})$.
Let $\tilde{P}_{r}$ be the polarization of $\opn{det}$.
One has
$\bsym{\Theta}^{1,1} = - \sum  \mrm{d} t_{i} \wedge \bsym{\theta}_{i}$
(cf.\ Lemma \ref{lem3.9}), so by Lemma \ref{lem3.2} we have
\begin{equation} \label{eqn7.3}
\opn{det} (\bsym{\Theta}^{1,1}) =
\tilde{P}_{r}(- \mrm{d} t_{0} \wedge \bsym{\theta}_{0}, \ldots,
- \mrm{d} t_{r-1} \wedge \bsym{\theta}_{r-1}) .
\end{equation}
But in case (iii), using Lemma \ref{lem2.2} again, we get
$\opn{det} (\bsym{\Theta}^{1,1}) = 0$.
\end{proof}

\begin{lem} \label{lem7.2}
Suppose $\xi = (x_{0}, \ldots, x_{r})$ is a saturated chain in $U$
satisfying: $x_{0}$ is the generic point of $X$, and
$a_{1}(x_{i}), \ldots, a_{i}(x_{i}) = 0$ for $0 \leq i \leq r$
\textup{(}so in particular $x_{r} \in Z$\textup{)}. Then
\[ \int_{\Delta^{r}} \opn{det} R_{\xi} =
(-1)^{\binom{r+1}{2}}
\frac{\mrm{d} a_{1} \wedge \cdots \wedge \mrm{d} a_{r}}{
a_{1} \cdots a_{r}} \in \Omega^{r}_{X / k, \xi} . \]
\end{lem}

\begin{proof}
Since the point $x_{r}$ falls into case 3 of (\ref{eqn7.2}), and for
every $i < r$, $x_{i}$  falls into case 2,
an easy linear algebra calculation shows that
$\bsym{\Theta}^{1, 1} = ( \mrm{d} t_{i-1} \wedge
a_{i}^{-1} \mrm{d} a_{j} )$
(i.e.\ $\mrm{d} t_{i-1} \wedge a_{i}^{-1} \mrm{d} a_{j}$ appears
in the $(i,j)$ position). By Lemma \ref{lem3.2},
\[ \begin{aligned}
\opn{det} \bsym{\Theta}^{1,1} & =
r! \mrm{d} t_{0} \wedge a_{1}^{-1} \mrm{d} a_{1}
\wedge \cdots \wedge \mrm{d} t_{r-1} \wedge
a_{r}^{-1} \mrm{d} a_{r} \\[1mm]
& = r! (-1)^{\binom{r}{2}}
\mrm{d} t_{0} \wedge \cdots \wedge \mrm{d} t_{r-1}
\wedge
\frac{\mrm{d} a_{1} \wedge \cdots \wedge \mrm{d} a_{r}}{
a_{1} \cdots a_{r}} .
\end{aligned} \]
Now use the fact that
$\int_{\Delta^{r}}
\mrm{d} t_{0} \wedge \cdots \wedge \mrm{d} t_{r-1} =
(-1)^{r} (r!)^{-1}$.
\end{proof}

\begin{lem} \label{lem7.3}
Let $\xi = (x_{0}, \ldots, x_{r} = z)$ be a saturated chain,
and let $\sigma : k(z) \ar \mcal{O}_{X, (z)}$ a coefficient field.
Then the Parshin residue map
$\opn{Res}_{k(\xi) / k(z)} :
\Omega^{\bdot, \mrm{sep}}_{k(\xi) / k} \ar
\Omega^{\bdot}_{k(z) / k}$
\textup{(}see \cite{Ye1} Definition \textup{4.1.3)} satisfies
\[ \opn{Res}_{k(\xi) / k(z)} ( \opn{dlog} a_{1} \wedge \cdots \wedge
\opn{dlog} a_{r} \wedge \alpha ) = 0 \]
for all $a_{1}, \ldots, a_{r} \in \mcal{O}_{X, (z)}$
and
$\alpha \in (\mfrak{m}_{z} + \mrm{d} \mfrak{m}_{z})
\Omega^{\bdot}_{X/k, (z)}$.
\end{lem}

\begin{proof}
By induction on $r$. We start with $r=1$.
By the definition of residues, it suffices to prove that for any
local factor $L$ of $k(\xi)$,
$\opn{Res}_{L / k(z)}(\opn{dlog} a_{1} \wedge \alpha) = 0$.
Now $L \cong K((t))$, with $K$ a finite field extension of $k(z)$,
and the image of $\mcal{O}_{X, (z)} \ar L$ lies in $K[\sqbr{t}]$
(we are using the fact that $\opn{char} k = 0$). Note that
$\Omega^{\bdot}_{X/k, (z)} =
\Omega^{\bdot, \mrm{sep}}_{\mcal{O}_{X, (z)} / k}$, so
$\alpha = t \beta + \mrm{d} t \wedge \gamma$ for some
$\beta, \gamma \in \Omega^{\bdot, \mrm{sep}}_{K[\sqbr{t}] / k}$.
Also, $a_{1} = t^{e} u$ with $u \in K[\sqbr{t}]^{*}$ and
$e \in \mbb{Z}$,
so $\opn{dlog} a_{1} = e \opn{dlog} t + \opn{dlog} u$. But
\[ \opn{Res}_{L / K} ((e \opn{dlog} t + \opn{dlog} u)
(t \beta + \mrm{d} t \wedge \gamma)) = 0 . \]

Now assume $r > 1$, and set
$\partial_{r} \xi := (x_{0}, \ldots, x_{r-1})$ and $y := x_{r-1}$.
First take $a_{1}, \ldots, a_{r}, \alpha$ algebraic, i.e.\
$a_{i} \in \mcal{O}_{X, z}$ and
$\alpha \in (\mfrak{m}_{z} + \mrm{d} \mfrak{m}_{z})
\Omega^{\bdot}_{X/k, z}$.
Let
$\tau: k(y) \ar \mcal{O}_{X, (y)}$ be a lifting compatible with
$\sigma : k(z) \ar \mcal{O}_{X, (z)}$ (cf.\ \cite{Ye1}
Definition 4.1.5; again we use $\opn{char} k = 0$), so by
\cite{Ye1} Corollary 4.1.16,
\[ \opn{Res}_{k(\xi) / k(z)} = \opn{Res}_{k((y,z)) / k(z)} \circ
\opn{Res}_{k(\partial_{r} \xi) / k(y)} :
\Omega^{\bdot}_{k(x) / k} \ar \Omega^{\bdot}_{k(z) / k} . \]
The lifting $\tau$ determines a decomposition
\[ \Omega^{\bdot, \mrm{sep}}_{\mcal{O}_{X, (y)} / k} =
\Omega^{\bdot}_{k(y) / k} \oplus (\mfrak{m}_{y} + \mrm{d}
\mfrak{m}_{y})
\Omega^{\bdot, \mrm{sep}}_{\mcal{O}_{X, (y)} / k} , \]
and we decompose $\alpha = \alpha_{0} + \alpha_{1}$
and
$\opn{dlog} a_{r} = \beta_{0} + \beta_{1}$
(or rather their images in
$\Omega^{\bdot, \mrm{sep}}_{\mcal{O}_{X, (y)} / k}$)
accordingly.
Using the $\Omega^{\bdot}_{k(y) / k}$-linearity of
$\opn{Res}_{k(\partial_{r} \xi) / k(y)}$ and induction applied to
$\beta_{0} \wedge \alpha_{1}$,
$\beta_{1} \wedge \alpha_{0}$ and
$\beta_{1} \wedge \alpha_{1}$, we get
\[ \opn{Res}_{k(\partial_{r} \xi) / k(y)} (
\opn{dlog} a_{1} \wedge \cdots \wedge
\opn{dlog} a_{r-1} \wedge (\opn{dlog} a_{r} \wedge \alpha)) =
m \beta_{0} \wedge \alpha_{0} \]
with $m \in \mbb{Z}$. Since $\beta_{0}$, $\alpha_{0}$ are respectively
the images
of $\opn{dlog} a_{r}$, $\alpha$ in
$\Omega^{\bdot, \mrm{sep}}_{k((y,z)) / k}$, again using induction
we have
\[ \opn{Res}_{k((y,z)) / k(z)} (\beta_{0} \wedge \alpha_{0}) = 0 . \]
Finally by the continuity of
$\opn{Res}_{k(\xi) / k(z)}$ the result holds for any
$a_{1}, \ldots, a_{r}, \alpha$.
\end{proof}

\begin{proof}[Proof of Theorem \ref{thm7.1}]
By the definition of the product and by Lemma \ref{lem7.1},
for evaluating the product
$C_{X} \cdot c_{r}(\mcal{E}, \nabla)$
we need only consider the components
$c_{r}(\mcal{E}, \nabla)_{\xi} = \int_{\Delta} \opn{det} R_{\xi}$
of $c_{r}(\mcal{E}, \nabla)$ for saturated chains
$\xi = (x_{0}, \ldots, x_{r})$, where $x_{0}$ is
the generic point of $X$, $x_{r} = z$ is the generic point of some
irreducible component $Z'$ of $Z$, and $a_{i}(x_{j}) = 0$ for all
$i \leq j$. Fix one such component of $Z'$, and
let $\Xi_{z}$ be the set of all such chains ending with $z$.
By definition of $C_{Z}$ we must show that the map
\begin{equation} \label{eqn7.4}
 (C_{X} \cdot \int_{\Delta} \opn{det} R)_{z} :
\Omega^{n-r}_{X / k, z} \ar \mcal{K}_{X}(z)
\end{equation}
factors through the maps
\[ \Omega^{n-r}_{X / k, z} \surj \Omega^{n-r}_{k(z) / k}
\xrightarrow{(-1)^{m} l} \Omega^{n-r}_{k(z) / k} =
\mcal{K}_{Z'_{\mrm{red}}}(z)
\subset \mcal{K}_{X}(z) , \]
where $Z'_{\mrm{red}}$ is the reduced scheme, $l$ is the length of the
artinian ring
$\mcal{O}_{Z', z} = \mcal{O}_{X, (z)} / (a_{1}, \ldots, a_{r})$, and
$m = \binom{r+1}{2} + nr$.

Choose a coefficient field $\sigma : k(z) \ar \mcal{O}_{X, (z)}$.
By Lemma \ref{lem7.2} and Example \ref{exa6.1} we have
for any $\alpha \in \Omega^{n-r}_{X / k, z}$:
\begin{equation} \label{eqn7.5}
\begin{split}
( C_{X} \cdot & c_{r}(\mcal{E}, \nabla) )_{z} (\alpha) \\[1mm]
& = (-1)^{m} \sum_{\xi \in \Xi_{z}}
( \mrm{d} a_{1} \wedge \cdots \wedge \mrm{d} a_{r}
\wedge \alpha)
\cdot \left(\frac{1}{a_{1} \cdots a_{r}} \right)_{\xi} \\[1mm]
& = (-1)^{m} \sum_{\xi \in \Xi_{z}}
\opn{Res}_{k(\xi) / k(z)} \left(
\frac{\mrm{d} a_{1} \wedge \cdots \wedge \mrm{d} a_{r} \wedge \alpha
}{a_{1} \cdots a_{r}} \right) .
\end{split}
\end{equation}
By Lemma \ref{lem7.3} we see that this expression vanishes for
$\alpha \in
\opn{Ker}(\Omega^{n-r}_{X / k, (z)} \surj \Omega^{n-r}_{k(z) / k})$,
so that we can assume
$\alpha  \in \opn{Im}(\sigma : \Omega^{n-r}_{k(z) / k} \ar
\Omega^{n-r}_{X / k, (z)})$. Now
$\opn{Res}_{k(\xi) / k(z)}$ is a graded left
$\Omega^{\bdot}_{k(z) / k}$-linear homomorphism of degree $-r$, so
$\alpha$ may be extracted. On the other hand, by
\cite{Hu2} Corollary 2.5 or \cite{SY} Theorem 0.2.5, and by
\cite{HK} Example 1.14.b, we get
\[ \sum_{\xi \in \Xi_{z}}
\opn{Res}_{k(\xi) / k(z)} \left(
\frac{\mrm{d} a_{1} \wedge \cdots \wedge \mrm{d} a_{r}
}{a_{1} \cdots a_{r}} \right) =
\opn{Res}_{\mcal{O}_{X, (z)} / k(z)}
\gfrac{ \mrm{d} a_{1} \wedge \cdots \wedge \mrm{d} a_{r} }{
a_{1}, \ldots, a_{r}}
= l . \]
This concludes the proof.
\end{proof}

\appendix
\section{Simplicial De Rham Theorem}

In this appendix we provide proofs to Theorems \ref{thm1.1}
and \ref{thm1.2}. These proofs are essentially contained in
\cite{BG} and \cite{HS1}, but not quite in the formulation
needed here. For notation see Section 1.

Suppose $M \in \Delta \msf{Mod}(\mbb{Q})$ and
$N \in \Delta^{\circ} \msf{Mod}(\mbb{Q})$.
For any $\sigma : [m] \ar [n]$ in $\Delta$ there are
homomorphisms
$\sigma^{*} : M^{m} \ar M^{n}$ and
$\sigma_{*} : N_{n} \ar N_{m}$.
Let
\begin{equation}
N \otimes_{\leftarrow} M \subset \prod_{n = 0}^{\infty}
\left( N_{n} \otimes_{\mbb{Q}} M^{n} \right)
\end{equation}
be the submodule consisting of all
$(u_{0}, u_{1}, \ldots)$, $u_{n} \in N_{n} \otimes_{\mbb{Q}} M^{n}$,
s.t.\ for each $\sigma : [m] \ar [n]$,
\[ (1 \otimes \sigma^{*})(u_{m}) = (\sigma_{*} \otimes 1)(u_{n}) . \]
Of course it suffices to check this condition for
$\sigma = \partial^{i}, s^{i}$.
If
$M \in \Delta \msf{DGMod}(\mbb{Q})$ and
$N \in \Delta^{\circ} \msf{DGMod}(\mbb{Q})$
we set
$(N \otimes_{\leftarrow} M)^{p,q} := N^{q,\bdot} \otimes_{\leftarrow}
M^{p,\bdot}$,
and in the usual way this yields
$N \otimes_{\leftarrow} M \in \msf{DGMod}(\mbb{Q})$.

\begin{exa} \label{exaA.1}
Taking $N = \Omega^{\bdot}(\Delta_{\mbb{Q}})$ we get
\[ \tilde{\mrm{N}} M =  \Omega^{\bdot}(\Delta_{\mbb{Q}})
\otimes_{\leftarrow} M . \]
\end{exa}

Given a simplicial set $S$, $\opn{Hom}_{\msf{Sets}}(S, \mbb{Q})$
is a cosimplicial $\mbb{Q}$-algebra. Let
\[ C^{\bdot}(S, \mbb{Q}) :=
\mrm{N} \opn{Hom}_{\msf{Sets}}(S, \mbb{Q})  \]
which is a DGA, called the algebra of normalized cochains on $S$.
Observe that
$C^{q}(S, \mbb{Q}) \cong
\opn{Hom}_{\msf{Sets}}(S_{q}^{\mrm{red}}, \mbb{Q})$,
where $S_{q}^{\mrm{red}}$ is the set of nondegenerate simplices.
In particular
$C^{\bdot}(\Delta^{n}, \mbb{Q})$ is a finite (noncommutative)
$\mbb{Q}$-algebra.
Define a DGA
\[ A^{\bdot}(S, \mbb{Q}) := \tilde{\mrm{N}}
\opn{Hom}_{\msf{Sets}}(S, \mbb{Q}) . \]
Applying this to the cosimplicial simplicial set
$\Delta = \{ \Delta^{n} \}$ we get simplicial DGAs
$C^{\bdot}(\Delta, \mbb{Q})$ and
$A^{\bdot}(\Delta, \mbb{Q})$.

\begin{lem} \label{lemA.1}
\begin{enumerate}
\item For any $M \in \Delta \msf{DGMod}(\mbb{Q})$ there an
isomorphism of DG modules
\[ \mrm{N} M \cong
C^{\bdot}(\Delta, \mbb{Q}) \otimes_{\leftarrow} M . \]
\item For any $S \in \Delta^{\circ} \msf{Sets}$ there are isomorphisms
of DGAs
\begin{eqnarray*}
A^{\bdot}(S, \mbb{Q}) & \cong &
\opn{Hom}_{\Delta^{\circ} \msf{Sets}}(S,
\Omega^{\bdot}(\Delta_{\mbb{Q}})) \\
C^{\bdot}(S, \mbb{Q}) & \cong &
\opn{Hom}_{\Delta^{\circ} \msf{Sets}}(S, C^{\bdot}(\Delta, \mbb{Q})) .
\end{eqnarray*}
\item There is an isomorphism of DGAs
\[ A^{\bdot}(\Delta^{n}, \mbb{Q}) \cong
\Omega^{\bdot}(\Delta^{n}_{\mbb{Q}}) . \]

\end{enumerate}
\end{lem}

\begin{proof}
1.\ Straightforward; cf.\ \cite{HS1} Prop.\ 3.1.3.

\noindent 2.\ For any finite $\mbb{Q}$-module $V$ one has
$\opn{Hom}_{\msf{Sets}}(S_{n}, V) \cong
\opn{Hom}_{\msf{Sets}}(S_{n}, \mbb{Q}) \otimes_{\mbb{Q}} V$.
Apply this to $V = \Omega^{q}(\Delta^{n}_{\mbb{Q}})$ and
$V = C^{q}(\Delta^{n}, \mbb{Q})$.

\noindent 3.\ This follows from part 2 for $S = \Delta^{n}$,
and using the isomorphism
$T_{n} \cong$ \linebreak
$\opn{Hom}_{\Delta^{\circ} \msf{Sets}}(\Delta^{n}, T)$
for any $T \in \Delta^{\circ} \msf{Sets}$.
\end{proof}

Regarding notation, the DGA  we call
$C^{\bdot}(S, \mbb{Q})$ (resp.\ $A^{\bdot}(S, \mbb{Q})$)
is denoted by $C^{*}(S)$ (resp.\ $A^{*}(S)$) in \cite{BG}.
Our simplicial
DGA $C^{\bdot}(\Delta, \mbb{Q})$ is denoted by $Z^{*}$ in
\cite{HS1}.

Define a homomorphism
\begin{equation} \label{eqnA.1}
\rho : \Omega^{\bdot}(\Delta^{n}_{\mbb{Q}}) \ar
C^{\bdot}(\Delta^{n}, \mbb{Q})
\end{equation}
in $\msf{DGMod}(\mbb{Q})$ by
$\rho (\alpha)(\sigma) = \int_{\Delta^{m}_{\mrm{top}}}
\sigma_{*}(\alpha) \in \mbb{Q}$
for a differential form
$\alpha \in \Omega^{\bdot}(\Delta^{n}_{\mbb{Q}})$
and a simplex
$\sigma : [m] \ar [n]$
(cf.\ \cite{BG} \S 2.1 or \cite{HS1} \S 4.4.1).
When $[n] \in \Delta$ varies $\rho$ becomes a transformation of
functors $\Delta^{\circ} \ar \msf{DGMod}(\mbb{Q})$, i.e.\ a morphism
in
$\Delta^{\circ} \msf{DGMod}(\mbb{Q})$.
One directly verifies:

\begin{lem} \label{lemA.2}
Let $M \in \Delta \msf{DGMod}(\mbb{Q})$. Under the identifications
$\tilde{\mrm{N}} M \cong
\Omega^{\bdot}(\Delta_{\mbb{Q}}) \otimes_{\leftarrow} M$
\textup{(}Example \textup{\ref{exaA.1})} and
$\mrm{N} M \cong C^{\bdot}(\Delta, \mbb{Q}) \otimes_{\leftarrow} M$
\textup{(}Lemma \textup{\ref{lemA.1})},
the homomorphism
$\int_{\Delta} : \tilde{\mrm{N}} M \ar \mrm{N} M$
corresponds to
$\rho \otimes_{\leftarrow} 1$.
\end{lem}

\begin{proof}[Proof of Theorem \ref{thm1.1}]
In the course of the proof of Theorem 2.2 of \cite{BG} it is shown,
using an acyclic models argument, that the homomorphism
\[ \rho : A^{\bdot}(S, \mbb{Q}) \ar C^{\bdot}(S, \mbb{Q}) \]
is a homotopy equivalence, of functors
$\Delta^{\circ} \msf{Sets} \ar \msf{DGMod}(\mbb{Q})$.
This means that there exist natural transformations
$\phi : C^{\bdot}(S, \mbb{Q}) \ar A^{\bdot}(S, \mbb{Q})$,
$h$ and $h'$ satisfying
$\rho \phi = h \mrm{D} + \mrm{D} h$ and
$\phi \rho = h' \mrm{D} + \mrm{D} h'$.
Taking $S = \Delta^{n}$ and using Lemma \ref{lemA.1} we conclude that
the morphism $\rho$ of equation (\ref{eqnA.1}) is a homotopy
equivalence in $\Delta^{\circ} \msf{DGMod}(\mbb{Q})$. Applying
$(-) \otimes_{\leftarrow} M$ it follows that
$\int_{\Delta} : \tilde{\mrm{N}} M \ar \mrm{N} M$
is a homotopy equivalence in $\msf{DGMod}(\mbb{Q})$.
\end{proof}

\begin{proof}[Proof of Theorem \ref{thm1.2}]
According to \cite{BG} Prop.\ 3.3 (and using our Lemma
\ref{lemA.1}), there exists a $\mbb{Q}$-linear homomorphism
\[ \rho_{2} : \Omega^{\bdot}(\Delta_{\mbb{Q}}) \otimes_{\mbb{Q}}
\Omega^{\bdot}(\Delta_{\mbb{Q}})  \ar
C^{\bdot}(\Delta, \mbb{Q}) \]
of functors on $\Delta^{\circ}$, of degree $-1$, s.t.\
\[ \mrm{d} \rho_{2} + \rho_{2} \mrm{d} =
\rho \mu - \mu (\rho \otimes \rho) . \]
Here $\mu$ is multiplication.
Extend this bilinearly to a graded homomorphism
\[ \rho_{2} : (\Omega^{\bdot}(\Delta_{\mbb{Q}}) \otimes_{\mbb{Q}} A)^{
\otimes 2} \ar C^{\bdot}(\Delta, \mbb{Q}) \otimes_{\mbb{Q}} A \]
i.e.\
\[ \rho_{2}((\alpha \otimes a) \otimes (\beta \otimes b)) =
(-1)^{p q'} \rho_{2}(\alpha \otimes \beta) \otimes a \cdot b  \]
for
$a \in A^{p,l}$, $b \in A^{p', l}$,
$\alpha \in \Omega^{q}(\Delta^{l}_{\mbb{Q}})$ and
$\beta \in \Omega^{q'}(\Delta^{l}_{\mbb{Q}})$.
Setting
\[ \rho_{2}(u \otimes v) := (\rho_{2}(u_{0} \otimes v_{0}),
\rho_{2}(u_{1} \otimes v_{1}), \ldots) \]
for
$u = (u_{0}, u_{1}, \ldots), v = (v_{0}, v_{1}, \ldots) \in
\tilde{\mrm{N}} A$
we get a $\mbb{Q}$-linear homomorphism
$\rho_{2} : (\tilde{\mrm{N}} A)^{\otimes 2} \ar \mrm{N} A$.
A simple calculation shows that for any pair of cocycles
$u, v \in \tilde{\mrm{N}} A$,
\[ \rho (u \otimes v) - \rho(u) \cdot \rho(v) =
\mrm{D} \rho_{2} (u \otimes v) . \]
\end{proof}



\end{document}